\documentclass[aos,MSNbibl,seceqn,dvips]{arximspdf}
\usepackage{mathrsfs}
\usepackage{multirow}
\usepackage{graphicx}

\doi{10.1214/12-AOS1006} 
\volume{40}
\issue{2}
\pubyear{2012}
\firstpage{1233}
\lastpage{1262}

\makeatletter
\renewcommand{\epsilon}{\varepsilon}

\newcommand{\Lvb}{L\'{e}vy\ }
\newcommand{\Lv}{L\'{e}vy}
\newcommand{\Ito}{It\^o\ }
\newcommand{\cadlag}{c\`{a}dl\`{a}g\ }

\newtheorem{ass}{Assumption}
\newproclaim{remark}{Remark}
\newtheorem{theorem}{Theorem}
\newproclaim{rem}[prop]{Remark}
\newproclaim{exam}[prop]{Example}
\newtheorem{lema}{Lemma}
\newcommand{\VXb}{VIX\ }

\makeatother

\begin{document}
\begin{frontmatter}

\title{Realized Laplace transforms for pure-jump semimartingales\thanksref{T1}}
\runtitle{Realized Laplace transforms for pure-jump}
\thankstext{T1}{Supported in part by NSF Grant SES-0957330.}

\begin{aug}
\author[A]{\fnms{Viktor} \snm{Todorov}\corref{}\ead[label=e1]{v-todorov@northwestern.edu}}
\and
\author[B]{\fnms{George} \snm{Tauchen}\ead[label=e2]{george.tauchen@duke.edu}}
\runauthor{V. Todorov and G. Tauchen}
\affiliation{Northwestern University and Duke University}
\address[A]{Department of Finance\\
Northwestern University\\
Evanston, Illinois 60208-2001\\
USA\\
\printead{e1}} 
\address[B]{Department of Economics\\
Duke University\\
Durham, North Carolina 27708-0097\\
USA\\
\printead{e2}}
\end{aug}

\received{\smonth{7} \syear{2011}}
\revised{\smonth{3} \syear{2012}}

%
\begin{abstract}
We consider specification and inference for the stochastic scale of
discretely-observed pure-jump
semimartingales with locally stable \Lvb densities in the setting where
both the time span of the data
set increases, and the mesh of the observation grid decreases. The
estimation is based on constructing
a nonparametric estimate for the empirical Laplace transform of the
stochastic scale over a given
interval of time by aggregating high-frequency increments of the
observed process on that time
interval into a statistic we call realized Laplace transform.
The realized Laplace transform depends on the activity of the driving
pure-jump martingale, and we
consider both cases when the latter is known or has to be inferred from
the data.
\end{abstract}

%

\end{frontmatter}
%

\section{Introduction}
Continuous-time semimartingales are used extensively for modeling many
processes in various areas, particularly in finance. Typically the
model of interest is an \Ito semimartingale (semimartingale with
absolute continuous characteristics) given by
%
\begin{equation}\label{eq:intro}
dX_t = \alpha_t\,dt+\sigma_{t-}\,dZ_t+dY_t,
\end{equation}
where $\alpha_t$ and $\sigma_t$ are some processes with \cadlag paths,
$Z_t$ is an infinite variation \Lvb martingale and $Y_t$ is a finite
variation jump process satisfying certain regularity conditions (all
technical conditions on the various processes will be given later). The
martingale $Z_t$ can be continuous (i.e., Brownian motion),
jump-diffusion or of pure-jump type (i.e., without a continuous
component). The presence of the last term in (\ref{eq:intro}) might
appear redundant, as $Z_t$ can already contain jumps, but its presence
will allow us to encompass also the class of time-changed \Lvb
processes in our analysis. In any case, this last term in
(\ref{eq:intro}) is dominated over small scales by the term involving
$Z_t$.\vadjust{\goodbreak}

Our interest in this paper will be in inference to the process
$\sigma_t$ when $Z_t$ is a~pure-jump \Lvb process. Pure-jump models
have been used to study various processes of interest such as
volatility and volume of financial prices~\cite{BNS01,Davis}, traffic
data~\cite{Mikosch} and electricity prices~\cite{Kluppelberg}.

Parametric or nonparametric estimation of a model satisfying
(\ref{eq:intro}) in the pure-jump case is quite complicated for at
least the following reasons. First, very often the transitional density
of $X_t$ is not known in closed-form. This holds true, even in the
relatively simple case when $X_t$ is a pure-jump \Lvb process. Second,
in many situations the realistic specification of $\sigma_t$ often
implies that $X_t$ is not a Markov process, with respect to its own
filtration, and hence all developed methods for estimation of the
latter will not apply. Third, the various parameters of the model
(\ref{eq:intro}) capture different statistical properties of the
process $X_t$ and hence will have various rates of convergence
depending on the sampling scheme. For example, in general, $\alpha_t$
and the tails of $Z_t$ can be estimated consistently only when the span
of the data increases, whereas the so-called activity of $Z_t$ can be
recovered even from a fixed-span data set, provided the mesh of the
latter decreases; see~\cite{B10} for estimation of the activity from
low-frequency data set. Finally, the simulation of the process $X_t$
can be, in many cases, difficult or time consuming.

In view of the above-mentioned difficulties, our goal here is
specification analysis for only part of the model, mainly the process
$\sigma_t$, in the case when~$X_t$ is a pure-jump process following
(\ref{eq:intro}). We conduct inference in the case when we have a
high-frequency data set of $X_t$ with increasing time span; see~\cite{B10} and~\cite{NR} for inference about jump processes based on
low frequency. We refer to $\sigma_t$ as the stochastic scale of the
pure-jump process $X_t$ in analogy with the scale parameter of a stable
process, that is, when $X_t$ is a stable process then the constant
$\sigma$ is the scale of the process. $\sigma_t$ is key in the
specification of (\ref{eq:intro}) and, in particular, it captures the
time-variation of the process $X_t$ over small intervals of time. Our
goal here will be to make the inference about $\sigma_t$ robust to the
rest of the components of the model, that is, the specification of
$\alpha_t$ and $Y_t$, as well as the dependence between $\sigma_t$ and~$Z_t$.

The inference in the paper is for processes for which the \Lvb measure
of the driving martingale $Z_t$ in (\ref{eq:intro}) behaves around zero
like that of a stable process. This covers, of course, the stable
process, but also many other \Lvb processes of interest with details
provided in Section~\ref{sec:setup} below. The idea of our proposed
method of inference is to use the fact that when $Z_t$ is locally
stable, the leading component of the process $X_t$ over small scales is
governed by that of the ``stable component'' of $Z_t$. Moreover, when
$\sigma_t$ is an \Ito semimartingale, then ``locally'' its changes are
negligible and $\sigma_t$ can be treated as constant. Intuitively then
infill asymptotics can be conducted as if the increments of $X_t$ are
products of (a locally constant) stochastic scale and independent
i.i.d. stable random variables. This, in particular, implies that the
empirical characteristic function of the high-frequency\vadjust{\goodbreak} increments over
a small interval of time will estimate the characteristic function of a~scaled stable process. The latter, however, is the Laplace transform
for the locally constant stochastic scale. Therefore, aggregating over
a fixed interval time the empirical characteristic function of the
(appropriately re-scaled) high-frequency increments of $X_t$ provides a
nonparametric estimate for the empirical Laplace transform of the
stochastic scale over that interval. We refer to this simple statistic
as realized Laplace transform for pure-jump processes. The connection
between the empirical characteristic function of the driving martingale
and the Laplace transform of the stochastic scale in the context of
time-changed \Lvb processes, with time-change independent of the
driving martingale, has been previously used for low-frequency
estimation in~\cite{B11}.

The inference based on the realized Laplace transform is robust to the
specification of $\alpha_t$, as well as the tail behavior of $Z_t$.
Intuitively, this is due to our use of the high-frequency data whose
marginal law is essentially determined by the small jumps of $Z_t$ and
the stochastic scale $\sigma_t$. Quite naturally, however, our
inference depends on the activity of the small jumps of the driving
martingale $Z_t$. The latter corresponds to the index of the stable
part of $Z$, and using the self-similarity of the stable process, it
determines its scaling over different (high) frequencies. Therefore,
the activity index enters directly into the calculation of the realized
Laplace transform. We conduct inference, both in the case where the
activity is assumed known, and when it needs to be estimated from the
data. The estimation of the activity index, however, differs from the
inference for the stochastic scale. While for the latter we need, in
general, the time-span to increase to infinity (except for the
degenerate case when $\sigma_t$ is actually constant), for the former,
this is not the case. The activity index can be estimated only with a
fixed span of high-frequency data, and, in general, increasing
time-span will not help for its nonparametric estimation. Therefore, we
estimate the activity index of~$Z_t$ using the initial part of the
sample with a fixed span and then plug it into the construction of the
realized Laplace transform. We further quantify the asymptotic effect
from this plug-in approach on the inference for the Laplace transform
of the stochastic scale.

The Laplace transform of the stochastic scale preserves the information
for its marginal distribution. Therefore, it can be used for efficient
estimation and specification testing. We illustrate this in a
parametric setting by minimizing the distance between our nonparametric
Laplace estimate and a model implied one, which is similar to
estimation based on the empirical characteristic function, as in
\cite{Feuerverger}.

Finally, the current paper studies the realized Laplace transform for
the case when $Z_t$ is pure-jump, while~\cite{TT10} (and the empirical
application of it in~\cite{TTG11}) consider the case where $Z_t$ is a
Brownian motion. The pure-jump case is substantively different,
starting from the very construction of the statistic as well as its\vadjust{\goodbreak}
asymptotic behavior. The leading component in the asymptotic expansions
in the pure-jump case is a stable process with an index of less than
$2$, and this index is, in general, unknown and needs to be estimated,
which further necessitates different statistical analysis from the
continuous case. Also, the residual components in $X_t$, like
$\alpha_t$, play a~more prominent role when $X_t$ is of pure-jump type,
and when the activity is low this requires modifying appropriately the
realized Laplace transform to purge them.

The paper is organized as follows. Section~\ref{sec:setup} presents the
formal setup and assumptions. Section~\ref{sec:estim} introduces the
realized Laplace transform and derives its limit behavior. In
Section~\ref{sec:mc} we conduct a Monte Carlo study, and in
Section~\ref{sec:param} we present a parametric application of the
developed limit theory. Section~\ref{sec:concl} concludes. The proofs
are given in Section~\ref{sec:proofs}.

\section{Setting and assumptions}\label{sec:setup}
Throughout the paper, the process of interest is denoted with $X_t$ and
is defined on some filtered probability space
$ (\Omega,\mathcal{F},(\mathcal{F}_t)_{t\geq0},
\mathbb{P} )$. Before stating our assumptions, we recall that a
\Lvb process $L_t$ with characteristic triplet $(b,c,\nu)$ with respect
to truncation function $\kappa(x) = x$ (we will always assume that the
process has a finite first order moment) is a process with
characteristic function given by
%
\begin{equation}
\mathbb{E} (e^{iuL_t} ) =
\exp \biggl(t \biggl(iub-u^2c/2+\int_{\mathbb{R}}(e^{iux}-1-iux)\nu
(dx) \biggr) \biggr).
\end{equation}
With this notation, we assume that the \Lvb process $Z_t$ in
(\ref{eq:intro}) has a~characteristic triplet $(0,0,\nu)$ for $\nu$
some \Lvb measure. Note that since the truncation function with respect
to which the characteristics of the \Lvb processes are presented is the
identity, the above implies that $Z_t$ is a pure-jump martingale. The
first term in (\ref{eq:intro}) is the drift term. It captures the
persistence in the process, and when $X_t$ is used to model financial
prices, the drift captures compensation for risk and time. The second
term in (\ref{eq:intro}) is defined in a stochastic sense, since in
Assumption~\ref{assa} below, we will assume that $Z_t$ is of infinite variation.
The last term in (\ref{eq:intro}) is a finite variation pure-jump
process. Assumption~\ref{assa} below will impose some restrictions on its
properties, but we stress that there is no assumption of independence
between the processes $\sigma_t$, $Z_t$ and $Y_t$.

In the pure-jump model the jump martingale, $Z_t$ substitutes the
Brownian motion used in jump-diffusions to model the ``small'' moves.
We note that the ``dominant'' part of the increment of $X_t$ over a
short interval of time $(t,t+\Delta)$ is $\sigma_{t-} \times
(Z_{t+\Delta} - Z_{t})$. This term is of order $O_p(\Delta^{\alpha})$
for $\alpha\in[1/2,1)$, while the rest of the components of $X_t$ are
at most $O_p(\Delta)$ when $\Delta\downarrow0$.

We recall from the introduction that our object of interest in this
paper is the stochastic scale of the martingale component of $X_t$,
that is, $\sigma_t$. Of course we observe only $X_t$ and $\sigma_t$ is
hidden into it, so our goal in the paper will be to uncover $\sigma_t$,
and its distribution in particular, with assuming as little as possible\vadjust{\goodbreak}
about the rest of the components of $X_t$ and the specification of
$\sigma_t$ itself (including the activity of the driving martingale).
Given the preceding discussion, the scaling of the driving martingale
components over short intervals of time will be of crucial importance
for us, as, at best, we can observe only a product of the stochastic
scale with $Z_{t,t+\Delta}$. Our Assumption~\ref{assa} below characterizes the
behavior of $Z_t$ and $Y_t$ over small scales.\looseness=-1

\def\theass{\Alph{ass}}
\begin{ass}\label{assa}
The \Lvb density of $Z_t$, $\nu$, is
given by
%
\begin{equation}\label{eq:ass_a_pj}
\nu(x) = \frac{A}{|x|^{\beta+1}}+\nu'(x),\qquad \beta\in
(1,2),\qquad
\int_{\mathbb{R}}|x|\nu(x)\,dx<\infty,
\end{equation}
where
%
\begin{equation}\label{eq:A}
A =  \biggl(\frac{4\Gamma(2-\beta)|\cos(\beta\pi
/2)|}{\beta(\beta-1)} \biggr)^{-1},\qquad \beta\in(1,2)
\end{equation}
and further there exists $x_0>0$ such that for $|x|\leq x_0$ we have
$|\nu'(x)|\leq\frac{C}{|x|^{\beta'+1}}$ for some $\beta'<1$ and a
constant $C\geq0$.

We further have $Y_t$ absolutely integrable and
$\mathbb{E}|Y_t-Y_s|^{\beta'} < C|t-s||\log|t-s||$ for every
$t, s\geq
0$ with $|t-s|\leq1$, some positive constant $C$, and $\beta'<1$ being
the constant above.
\end{ass}

Assumption~\ref{assa} implies that the small scale behavior of the driving
martingale $Z_t$ is like that of a stable process with index $\beta$.
The index $\beta$ determines the ``activity'' of the driving process,
that is, the vibrancy of its trajectories, and thus henceforth we will
refer to it as the activity. Formally, $\beta$ equals the
Blumenthal--Getoor index of the \Lvb process $Z_t$. The value of the
index~$\beta$ is crucial for recovering $\sigma_t$ from the discrete
data on $X_t$, as intuitively it determines how big on average the
increments $Z_{t,t+\Delta}$ should be for a given sampling frequency.
The following lemma makes this formal.

\begin{lema}\label{lema_small_scale}
Let $Z_t$ satisfy Assumption~\ref{assa}. Then for $h\rightarrow0$, we have
%
\begin{equation}\label{eq:small_scale_1}
h^{-1/\beta}Z_{th} \stackrel{\mathcal{L}}{\longrightarrow} S_t,
\end{equation}
where the convergence is for the Skorokhod topology on the space of
\cadlag functions. $\mathbb{R}_+\rightarrow\mathbb{R}$, and $S_t$ is a
stable process with characteristic function
$\mathbb{E} (e^{iuS_1} )=e^{-|u|^{\beta}/2}$.
\end{lema}

The value of the constant $A$ in (\ref{eq:A}) is a normalization that
we impose. We are obviously free to do that since what we observe is
$X_t$, whose leading component over small scales is an integral of
$\sigma_t$, with respect to the jump martingale $Z_t$, and we never
observe the two separately. The above choice of $A$ is a convenient one
that ensures that when $\beta\rightarrow2$, the jump process converges
finite-dimensionally to Brownian motion. We note that in Assumption~\ref{assa}
we rule out the case $\beta\leq1$, but this is done for brevity of
exposition, as most processes of interest are of infinite variation
(although we rule out some important processes like the generalized
hyperbolic).\vadjust{\goodbreak}

In Assumption~\ref{assa} we restrict the ``activity'' of the ``residual'' jump
components of $X_t$; that is, we limit their effect in determining the
small moves of $X$. The effect of the ``residual'' jump components on
the small moves is controlled by the parameter $\beta'$. From
(\ref{eq:ass_a_pj}), the leading component of $\nu(x)$ is the \Lvb
density of a stable, and $\nu'(x)$ is the residual one. The restriction
$\beta'<1$, implies that the ``residual'' jump component is of finite
variation. This restriction is not necessary for convergence in
probability results (only $\beta'<\beta$ is needed for this), but is
probably unavoidable if one needs also the asymptotic distribution of
the statistics that we introduce in the paper. In most parametric
models this restriction is satisfied.

We note that $\nu'(x)$ in (\ref{eq:ass_a_pj}) is a signed measure, and
therefore Assumption~\ref{assa} restricts only the behavior of $\nu(x)$ for
$x\sim0$ to be like that of a stable process. However, for the big
jumps, that is, when $|x|>K$ for some arbitrary $K>0$, the stable part
of $\nu(x)$ can be completely eliminated or tempered by negative values
of the ``residual'' $\nu'(x)$. An example of this, which is covered by
our Assumption~\ref{assa}, is the tempered stable process of~\cite{Ro04},
generated from the stable by tempering its tails, which has all of its
moments finite. Therefore, while Assumption~\ref{assa} ties the small scale
behavior of the driving martingale $Z_t$ with that of a stable process,
it leaves its large scale behavior unrestricted (i.e., the limit of
$h^{-\alpha}Z$th for some $\alpha>0$ when $h\rightarrow\infty$ is
unrestricted by our assumption) and thus, in particular, unrelated with
that of a stable process.\vspace*{-3pt}

\begin{remark}
Assumption~\ref{assa} is analogous to the
assumption used in~\cite{SJ07}. It is also related to the so-called
regular \Lvb processes of exponential type, studied in~\cite{BL02},
with $\beta=\nu$ in the notation of that paper. Compared with the above
mentioned processes of~\cite{BL02}, we impose slightly more structure
on the \Lvb density around zero but no restriction outside of it. We
note that if Assumption~\ref{assa} fails, then the results that follow are not
true. The degree of the violation depends on the sampling frequency and
the deviation of the characteristic function of $Z$ over small scales
from that of a stable.

Finally, the process $Y_t$ also captures a ``residual'' jump component
of $X_t$ in terms of its small scale behavior. Assumption~\ref{assa} limits its
activity by~$\beta'$. The component of $Z$ corresponding to $\nu'(x)$
and $Y_t$ control the jump measure of~$X_t$ away from zero. Unlike the
former, whose time variation is determined by $\sigma_t$, the latter
has essentially unrestricted time variation. There is clearly some
``redundancy'' in the specification in (\ref{eq:intro}) in terms of
modeling the jumps of $Z_t$ away from zero, but this is done to cover
more general pure-jump models, as we make clear from the following two
remarks.\vspace*{-3pt}
\end{remark}

\begin{remark}\label{rem2}
Assumption~\ref{assa} nests time-changed \Lvb
processes with absolute continuous time changes (see, e.g.,
\cite{CGMY03}), that is, specifications of the form\vspace*{-3pt}
%
\begin{equation}\label{eq:tc}
dX_t = \alpha_t\,dt+\int_0^t\int_{\mathbb{R}}x\widetilde{\mu}(ds,dx),\vadjust{\goodbreak}
\end{equation}
where $\mu$ is a random integer-valued measure with compensator
$a_t\,dt\otimes\nu(dx)$ for some nonnegative process $a_t$ and \Lvb
measure $\nu(dx)$ satisfying (\ref{eq:ass_a_pj}) of Assumption~\ref{assa} and
$\widetilde{\mu} = \mu- \nu$. This can be shown using Theorem 14.68 of
\cite{Jacod79} or Theorem 2.1.2 of~\cite{JP}, linking integrals of
random functions with respect to Poisson measure and random
integer-valued measures, and implies that~$X_t$ in (\ref{eq:tc}) can be
equivalently represented as (\ref{eq:intro}) with $\sigma_t$ given by~$a_t^{1/\beta}$.

On the other hand, if we start with $X$ given by (\ref{eq:intro}), with
$Z$ strictly stable and no $Y$, we can show, using the definition of
jump compensator and Theorem II.1.8 of~\cite{JS}, that the latter is a
time-changed \Lvb process with time change $|\sigma_{t-}|^{\beta}$. For
more general ``stable-like'' \Lvb processes, we need to introduce an
additional term [this is $Y_t$ in (\ref{eq:intro})] in addition to the
above time-changed stable process.

We note that the connection of (\ref{eq:intro}) with the time-changed
\Lvb processes does not depend on the presence of any dependence
between $\sigma_t$ and $Z_t$.
\end{remark}

\begin{remark}
Assumption~\ref{assa} is also satisfied by the
pure-jump \Lv-driven CARMA models (continuous-time autoregressive
moving average) which have been used for modeling series exhibiting
persistence; see, for example,~\cite{BR01a} and the many references
therein. For these processes $\sigma_t$ in (\ref{eq:intro}) is a
constant.

Our next assumption imposes minimal integrability conditions on~$\alpha_t$ and~$\sigma_t$ and further limits the amount of variation in
these processes over short periods of time. Intuitively, we will need
the latter to guarantee that by sampling frequently enough we can treat
``locally'' $\sigma_t$ (and $\alpha_t$) as constant.
\end{remark}

\begin{ass}\label{assb}
The process $\sigma_t$ is an \Ito
semimartingale given by
%
\begin{equation}\label{eq:ass_b}
\sigma_t=\sigma_0+\int_0^t\tilde{\alpha}_s\,ds+\int_0^t\tilde
{\sigma}_s\,dW_s+\int_0^t\int_{\mathbb{R}}\underline{\delta}(s-,x)
\underline{\tilde{\mu}}(ds,dx),
\end{equation}
where $W$ is a Brownian motion; $\underline{\mu}$ is a homogenous
Poisson measure, with \Lvb measure $\underline{\nu}(dx)$, having
arbitrary dependence with $\mu$, for $\mu$ being the jump measure of
$Z$. We assume $|\underline{\delta}(t,x)|\leq
\gamma_t\underline{\delta}(x)$ for some integrable process~$\gamma_t$
and $\int_{\mathbb{R}}(\underline{\delta}(x)\wedge
1)^{\beta^{\prime\prime}}\underline{\nu}(dx)<\infty$ for some $\beta^{\prime\prime}<2$.
Further, for every $t$ and $s$, we have
%
\begin{equation}\label{eq:ass_b_1}
\cases{\displaystyle
\mathbb{E} \biggl(\alpha_t^{2}+\sigma_t^{2}+\tilde
{\alpha}_t^{2}+\tilde{\sigma}_t^2+(\tilde{\sigma}_t')^2+\int
_{\mathbb{R}}\underline{\delta}^2(t,x)
\underline{\nu}(dx) \biggr)<C,\cr
\displaystyle\mathbb{E} \biggl(|\alpha_t-\alpha_s|^2+|\tilde{\sigma}_t-\tilde
{\sigma}_s|^2+\int_{\mathbb{R}}\bigl(\underline{\delta}(t,x)-\underline
{\delta}(s,x)\bigr)^2\underline{\nu}(dx) \biggr)
<C|t-s|,
}\hspace*{-35pt}
\end{equation}
where $C>0$ is some constant that does not depend on $t$ and $s$.
\end{ass}

Assumption~\ref{assb} imposes $\sigma_t$ to be an \Ito semimartingale. This is a
relatively mild assumption satisfied by the popular multifactor affine
jump-diffusions~\cite{DFS03} as well as the CARMA \Lv-driven models
used to model persistent processes~\cite{BR01a}. Assumption~\ref{assb} rules
out certain long-memory specifications for $\sigma_t$, although we
believe that, at least for some of them, the results in this paper will
continue to hold.

Importantly, however, Assumption~\ref{assb} allows for jumps in the stochastic
scale that can have arbitrary dependence with the jumps in $X_t$ which
is particularly relevant for modeling financial data, for example, in
the parametric models of~\cite{COGARCH}. Finally, the second part of
(\ref{eq:ass_b_1}) will be satisfied when the corresponding processes
are \Ito semimartingales. The next Assumption~\ref{assb1} restricts Assumption~\ref{assb}
in a way that will allow us to strengthen some of the theoretical
results in the next section.

{\setcounter{ass}{1}
\def\theass{\Alph{ass}$^{\prime}$}
\begin{ass}\label{assb1}
The process $\sigma_t$ is an \Ito
semimartingale given by
%
%
\begin{equation}\label{eq:ass_b'}
\sigma_t=\sigma_0+\int_0^t\!\tilde{\alpha}_s\,ds+\int_0^t\!\tilde
{\sigma}_s\,dW_s+\int_0^t\!\tilde{\sigma}'_s\,dZ_s+\int_0^t\int
_{\mathbb{R}}\underline{\delta}(s-,x)
\underline{\tilde{\mu}}(ds,dx),\hspace*{-35pt}
\end{equation}
with the same notation as Assumption~\ref{assb}, with the only difference being
that~$\underline{\mu}$ is now independent from~$\mu$. We also assume
that corresponding condition~(\ref{eq:ass_b_1}) holds as well as that
$Z_t$ is square-integrable.
\end{ass}}

The strengthening in Assumption~\ref{assb1} is in the modeling of the
dependence between the jumps in $X_t$ and $\sigma_t$. In Assumption~\ref{assb1}, this is done via the third integral in (\ref{eq:ass_b'}). This is
similar to modeling dependence between continuous martingales using
correlated Brownian motions. What is ruled out by Assumption~\ref{assb1} is
dependence between the jumps in $X_t$ and $\sigma_t$ that is different
for the jumps of different size. Assumption~\ref{assb1} will be satisfied when
the pair $(X_t,\sigma_t)$ are modeled jointly via a \Lv-driven SDE.

Finally, in our estimation we make use of long-span asymptotics for the
process $\sigma_t$ and the latter contains temporal dependence.
Therefore, we need a condition on this dependence that guarantees that
a central limit theorem for the associated empirical process exists.
This condition is given next.

\begin{ass}\label{assc}
The volatility $\sigma_t$ is a
stationary and $\alpha$-mixing process with
$\alpha_t^{\mathrm{mix}}=O(t^{-3-\iota})$ for arbitrarily small
$\iota>0$ when $t\rightarrow\infty$, where
\begin{eqnarray*}
\alpha_t^{\mathrm{mix}}&=&\sup_{A\in\mathscr{F}_0, B\in\mathscr
{F}^t}|\mathbb{P}(A\cap
B)-\mathbb{P}(A)\mathbb{P}(B)|,\\  \mathscr{F}_0&=&\sigma(\sigma_s,
s\leq0) \qquad\mbox{and}\qquad
  \mathscr{F}^t=\sigma(\sigma_s, s\geq
t).
\end{eqnarray*}
\end{ass}

\section{Limit theory for RLT of pure-jump semimartingales}\label{sec:estim}
Now we are ready to formally define the realized Laplace transform for
the pure-jump model and derive its asymptotic properties. We assume the
process $X_t$ is observed at the equidistant times
$0,\Delta_n,\ldots,i\Delta_n,\ldots,[T/\Delta_n]$ where $\Delta_n$ is the
length of the high-frequency interval, and $T$ is the span of the data.\vadjust{\goodbreak}
The realized Laplace transform is then defined as
%
\begin{eqnarray}\label{eq:V}
V_T(X,\Delta_n,\beta,u)=\sum_{i=1}^{[T/\Delta_n]}\Delta_n\cos
((2u)^{1/\beta}\Delta_n^{-1/\beta}\Delta_i^nX),\nonumber \\[-8pt]\\ [-8pt]
\eqntext{\Delta_i^nX=X_{i\Delta_n}-X_{(i-1)\Delta_n},}
\end{eqnarray}
where $\beta$ is the activity index of the driving martingale $Z_t$
given in its \Lvb density in (\ref{eq:ass_a_pj}).
$V_T(X,\Delta_n,\beta,u)$ is the real part of the empirical
characteristic function of the appropriately scaled increments of the
process $X_t$. In the case of jump-diffusions, $\beta$ in (\ref{eq:V})
is replaced with $2$ as the activity of the Brownian motion has an
index of $2$ (i.e., for the Brownian motion,
Lemma~\ref{lema_small_scale} holds with $\beta$ replaced by $2$). We
show in this section that $V_T(X,\Delta_n,\beta,u)/T$ is a consistent
estimator for the empirical Laplace transform of $|\sigma_t|^{\beta}$
and further derive its asymptotic properties under various sampling
schemes as well as assumptions regarding whether $\beta$ is known or
needs to be estimated.\looseness=1

\subsection{Fixed-span asymptotics}
We start with the case when $T$ is fixed and $\Delta_n\rightarrow0$,
that is, the infill asymptotics, and we further assume we know $\beta$.
Since the driving martingale over small scales behaves like
$\beta$-stable (Assumption~\ref{assa}) and the stochastic scale changes over
short intervals are not too big on average (Assumption~\ref{assb}), then the
``dominant'' part (in a infill asymptotic sense) of the increment
$\Delta_i^nX$ (when $\Delta_n$ is small) is
$\sigma_{(i-1)\Delta_n}\Delta_i^nZ$. $\Delta_i^nZ$ is approximately
stable, and from Lemma~\ref{lema_small_scale}, we have approximately
$\Delta_i^nZ \stackrel{d}{=} \Delta_n^{1/\beta}\times Z_{1}$ with the
characteristic function of $Z_1$ given by $e^{-|u|^{\beta}/2}$.
Therefore, for a fixed $T$, $V_T(X,\Delta_n,\beta,u)$ is approximately
a sample average of a heteroscedastic data series. Thus, by the law of
large numbers (when $\Delta_n\rightarrow0$), it will converge to
$\int_0^Te^{-u|\sigma_t|^{\beta}}\,ds$, which is the empirical Laplace
transform of $|\sigma_t|^{\beta}$ after dividing by $T$. The following
theorem gives the precise infill asymptotic result. In it we denote
with $\mathcal{L}-s$ convergence stable in law, which means that the
convergence in law holds jointly with any random variable defined on
the original probability space.

\begin{theorem}\label{theorem:fixed}
For the process $X_t$, assume that Assumptions \textup{\ref{assa}} and \textup{\ref{assb}} hold with
$\beta'<\beta/2$, and let $\Delta_n\rightarrow0$ with $T$ fixed.
\begin{longlist}
\item[(a)] If $\beta>4/3$, then we have
%
\begin{eqnarray}\label{eq:fixed_1}
&&\frac{1}{\sqrt{\Delta_n}} \biggl(V_T(X,\Delta_n,\beta,u)-\int
_0^Te^{-u|\sigma_s|^{\beta}}\,ds \biggr)\nonumber\\ [-8pt]\\ [-8pt]
&&\qquad \stackrel{\mathcal
{L}-s}{\longrightarrow}
\sqrt{\int_0^TF_{\beta}(u^{1/\beta}|\sigma_s|)\,ds}\times E,\nonumber
\end{eqnarray}
where $E$ is a standard normal variable defined on an extension of the
original probability space and independent from the $\sigma$-field
$\mathcal{F}$;
$F_{\beta}(x)=\frac{e^{-2^{\beta-1}x^{\beta}}-2e^{-x^{\beta}}+1}{2}$
for $x>0$.

A consistent estimator for the asymptotic variance is given by
%
\begin{equation}\label{eq:fixed_2}
\frac{V_T(X,\Delta_n,\beta,2^{\beta-1}u)-2V_T(X,\Delta_n,\beta,u)+1}{2}.
\end{equation}
\item[(b)] If $\beta\leq4/3$, then
%
\begin{equation}\label{eq:fixed_3}
 \biggl(V_T(X,\Delta_n,\beta,u)-\int_0^Te^{-u|\sigma_s|^{\beta
}}\,ds \biggr)=O_p (\Delta_n^{2-2/\beta} ).
\end{equation}
\end{longlist}
\end{theorem}

In the case when $X_t$ is a \Lvb process, or, more generally, when only
the scale $\sigma_t$ is constant as is the case for the \Lv-driven
CARMA models, the above theorem can be used to estimate the scale
coefficient $\sigma$ by either fixing some $u$ or using a whole range
of $u$'s as in the methods for estimation of stable processes based on
the empirical characteristic function; see, for example,
\cite{Koutrouvelis} and~\cite{SJ06}. Furthermore, this can be done
jointly with the nonparametric estimation of $\beta$ by using for
example the estimator we proposed in~\cite{TT09} that we define in
(\ref{eq:bg_est}) below.

The limit result in Theorem~\ref{theorem:fixed} is driven by the small
jumps in $X_t$, and this allows us to disentangle the stochastic scale
(which drives their temporal variation) from the other components of
the model, mainly the jumps away from zero. This is due to the fact
that the cosine function is bounded and infinitely differentiable which
limits the effect of the jumps of size away from zero on it. By
contrast, for example, the infill asymptotic limit of the quadratic
variation of the discretized process is the quadratic variation of~$X_t$ which is determined by all jumps, not just the infinitely small
ones.

Unfortunately when the activity of the driving martingale is relatively
low, i.e., $\beta<4/3$, we do not have a CLT for
$V_T(X,\Delta_n,\beta,u)$. The reason is in the presence of the drift
term, which for the purposes of our estimation starts behaving closer
to the driving martingale $Z_t$ and this slows the rate of convergence
of our statistic. However, we can use the fact that over successive
short intervals of time the contribution of the drift term in the
increments of $X_t$ is the same while sum or difference of i.i.d.
stable random variables continues to have a stable distribution.
Therefore, if we difference the increments of $X_t$, we will remove the
drift term (up to the effect due to the time-variation in it which will
be negligible) and the leading term will still be a product of the
locally constant stochastic scale and a stable variable. Thus we
consider the following alternative estimator:
%
\begin{equation}\label{eq:L}\qquad
\widetilde{V}_T(X,\Delta_n,\beta,u)=\sum_{i=2}^{[T/\Delta
_n]}\Delta_n\cos [u^{1/\beta}\Delta_n^{-1/\beta}(\Delta
_i^nX-\Delta_{i-1}^nX) ].
\end{equation}

\begin{theorem}\label{theorem:fixed-2}
For the process $X_t$, assume that Assumptions \textup{\ref{assa}} and \textup{\ref{assb}} hold with
$\beta'<\beta/2$, and let $\Delta_n\rightarrow0$ with $T$ fixed. We
have
%
\begin{eqnarray}
&&\frac{1}{\sqrt{\Delta_n}} \biggl(\widetilde{V}_T(X,\Delta_n,\beta
,u)-\int_0^Te^{-u|\sigma_{s}|^{\beta}}\,ds \biggr)\nonumber\\ [-8pt]\\ [-8pt]
&&\qquad \stackrel{\mathcal
{L}-s}{\longrightarrow}
\sqrt{\int_0^T\widetilde{F}_{\beta}(u^{1/\beta}|\sigma
_s|)\,ds}\times
\widetilde{E},\nonumber
\end{eqnarray}
where $\widetilde{E}$ is a standard normal variable defined on an
extension of the original probability space and independent from the
$\sigma$-field $\mathcal{F}$; $\widetilde{F}_{\beta}(x) =
 (\frac{e^{-2^{\beta-1}x^{\beta}}+1}{2} )^2+2e^{-x^{\beta
}}\frac{e^{-2^{\beta-1}x^{\beta}}+1}{2}-3e^{-2x^{\beta}}+
 (\frac{e^{-2^{\beta-1}x^{\beta}}-1}{2} )^2$ for $x>0$.
\end{theorem}

A comparison of the standard errors in (\ref{theorem:fixed}) and
(\ref{theorem:fixed-2}) shows that the latter can be up to $2.5$ times
higher than the former for values of $\beta>4/3$. This is the cost of
removing the effect of the drift term via the differencing of the
increments. Therefore, $\widetilde{V}_T(X,\Delta_n,\beta,u)$ should be
used only in the case when $\beta\leq4/3$. For brevity, the results
that follow will be presented only for $V_T(X,\Delta_n,\beta,u)$, but
analogous results will hold for $\widetilde{V}_T(X,\Delta_n,\beta,u)$.

\subsection{Long span asymptotics: The case of known activity}
We continue next with the case when the time span of the data increases
together with the mesh of observation grid decreasing. The
high-frequency data allows us to ``integrate out'' the increments
$\Delta_i^nZ$, that is, it essentially allows to ``deconvolute''
$\sigma_t$ from the driving martingale of $X_t$ in a robust way. After
dividing by $T$, the infill asymptotic limit of (\ref{eq:V}) is the
empirical Laplace transform of the stochastic\vspace*{1pt} scale and we henceforth
denote it as $\widehat{\mathcal{L}}_{\beta}(u) =
\frac{1}{T}V_T(X,\Delta_n,\beta,u)$. Then, by letting
$T\rightarrow\infty$ we can eliminate the sampling variation due to the
stochastic nature of $\sigma_t$, and thus recover its population
properties, that is, estimate
$\mathcal{L}_{\beta}(u)=\mathbb{E} (e^{-u|\sigma_t|^{\beta
}} )$
which is the Laplace transform of $|\sigma_t|^{\beta}$.

The next theorem gives the asymptotic behavior of
$\widehat{\mathcal{L}}_{\beta}(u)$ when both $T\rightarrow\infty$ and
$\Delta_n\rightarrow0$. To state the result we first introduce some
more notation. We henceforth use the shorthand
%
\begin{equation}\label{eq:z_hat}\qquad
\widehat{Z}_{t,\beta}(u)=V_t(X,\Delta_n,\beta,u)-V_{t-1}(X,\Delta
_n,\beta,u),\qquad
t=1,\ldots,T,
\end{equation}
for the RLT over the time interval $(t-1,t)$. We further set\vspace*{-1pt}
%
%
\begin{equation}\label{eq:c_hat}
\widehat{C}_{k,\beta}(u,v)=\frac{1}{T}\sum_{t=k+1}^T \bigl(\widehat
{Z}_{t,\beta}(u)-\widehat{\mathcal{L}}_{\beta}(u) \bigr)
 \bigl(\widehat{Z}_{t-k,\beta}(v)-\widehat{\mathcal{L}}_{\beta
}(v) \bigr),\qquad  k\in\mathbb{N}.\hspace*{-35pt}
\end{equation}

\begin{theorem}\label{theorem:clt}
Suppose $T\rightarrow\infty$ and $\Delta_n\rightarrow0$, and the
process $X_t$ satisfies Assumptions \textup{\ref{assa}}, \textup{\ref{assb}} and \textup{\ref{assc}}.
\begin{longlist}
\item[(a)] We have
%
\begin{equation}\label{eq:clt_1}\qquad
\sqrt{T} \bigl(\widehat{\mathcal{L}}_{\beta}(u)-\mathcal{L}_{\beta
}(u) \bigr)=Y_T^{(1)}(u)+Y_T^{(2)}(u),
\end{equation}\vspace*{-10pt}
%
\begin{equation}\label{eq:clt_2}\qquad
 \cases{
Y_T^{(1)}(u) \stackrel{\mathcal{L}}{\longrightarrow} \Psi(u), \cr
Y_T^{(2)}(u)=O_p \bigl(\sqrt{T} \bigl(|\log\Delta_n|\Delta
_n^{1-\beta'/\beta}\vee\Delta_n^{
(2-2/\beta)\wedge
1/2} \bigr)\vee\sqrt{\Delta_n} \bigr),
}
\end{equation}
where the result for $Y_T^{(2)}(u)$ holds locally uniformly in $u$, and
the convergence of $Y_T^{(1)}(u)$ is on the space
$\mathcal{C}(\mathbb{R}_+)$ of continuous functions indexed by $u$ and
equipped with the local uniform topology (i.e., uniformly over compact
sets of $u\in\mathbb{R}_+$), and $\Psi(u)$ is a Gaussian process with
variance--covariance for $u, v>0$ given by
%
\begin{eqnarray}\label{eq:uniform_1}
\Sigma_{\beta}(u,v) &=&
\int_0^{\infty}\mathbb{E} \bigl[ \bigl(e^{-u|\sigma_t|^{\beta
}}-\mathcal{L}_{\beta}(u) \bigr)
 \bigl(e^{-v|\sigma_0|^{\beta}}-\mathcal{L}_{\beta}(v)
 \bigr)\nonumber\\ [-8pt]\\ [-8pt]
&&\phantom{\int_0^{\infty}\mathbb{E} \bigl[}{}+ \bigl(e^{-v|\sigma_t|^{\beta}}-\mathcal{L}_{\beta
}(v) \bigr)
 \bigl(e^{-u|\sigma_0|^{\beta}}-\mathcal{L}_{\beta}(u) \bigr) \bigr]\,dt.\nonumber
\end{eqnarray}
If we strengthen Assumption~\ref{assb} to Assumption~\ref{assb1}, we get the stronger
\[
Y_T^{(2)}(u)=O_p \bigl(\sqrt{T} \bigl(|\log\Delta_n|\Delta
_n^{1-\beta'/\beta}\vee\Delta_n^{2-2/\beta} \bigr)\vee\sqrt
{\Delta_n} \bigr).
\]
\item[(b)] For arbitrary integer $k\geq1$ and every $u,v>0$ we have
%
\begin{equation}\label{eq:se_0}
\widehat{Z}_{1,\beta}(u)\widehat{Z}_{k,\beta}(v) \stackrel{\mathbb
{P}}{\longrightarrow} \int_0^1\int_{k-1}^{k}
e^{-u|\sigma_t|^{\beta}}e^{-v|\sigma_s|^{\beta}}\,ds\,dt.
\end{equation}
If, further, $L_T$ is a deterministic sequence of integers satisfying
$\frac{L_T}{\sqrt{T}}\rightarrow0$ as $T\rightarrow\infty$ and
$L_T (|\log\Delta_n|\Delta_n^{1-\beta'/\beta}\vee\Delta_n^{
(2-2/\beta)\wedge1/2} )\rightarrow0$, we have
%
\begin{eqnarray}\label{eq:se_1}
\widehat{\Sigma}_{\beta}(u,v)&=&\widehat{C}_{0,\beta}(u,v)+\sum
_{i=1}^{L_T}\omega(i,L_T)\bigl(\widehat{C}_{i,\beta}(u,v)+\widehat
{C}_{i,\beta}(v,u)\bigr)\nonumber\\ [-8pt]\\ [-8pt]
&\stackrel{\mathbb{P}}{\longrightarrow}& \Sigma_{\beta}(u,v),\nonumber
\end{eqnarray}
where $\omega(i,L_T)$ is either a Bartlett or a Parzen kernel.
\end{longlist}
\end{theorem}

The result in (\ref{eq:clt_1}) holds locally uniformly in $u$. This is
important as, in a~typical application, one needs the Laplace transform
as a function of~$u$. We illustrate, in the next section, an
application of the above result to parametric estimation that makes use
of the uniformity. We note also that $\Sigma_{\beta}(u,v)$ is well
defined because of Assumption~\ref{assc}; see~\cite{JS}, Theorem VIII.3.79.

Under the conditions of Theorem~\ref{theorem:clt}, the scaled and
centered realized Laplace transform can be split into two components,
$Y_T^{(1)}(u)$ and $Y_T^{(2)}(u)$, that have different asymptotic
behavior and capture different errors involved in the estimation. The
first one, $Y_T^{(1)}(u)$, equals
$\sqrt{T} (\frac{1}{T}\int_0^Te^{-u|\sigma_t|^{\beta
}}\,dt-\mathbb{E} (e^{-u|\sigma_t|^{\beta}} ) )$,
which is the empirical process corresponding to the case of
continuous-record of $X_t$ in which case $|\sigma_t|^{\beta}$ can be
recovered exactly. Hence the magnitude of $Y_T^{(1)}(u)$ is the sole
function of the time span $T$. On the other hand, the term
$Y_T^{(2)}(u)$ captures the effect from the discretization error, that
is, the fact that we use high-frequency data and not continuous record
of $X_t$ in the estimation. For $Y_T^{(2)}(u)$ to be negligible, we
need a condition for the relative speed of $\Delta_n\rightarrow0$ and
$T\rightarrow\infty$ which, in the general case of Assumption~\ref{assb}, is
given by
$\sqrt{T} (|\log\Delta_n|\Delta_n^{1-\beta'/\beta}\vee
\Delta_n^{
(2-2/\beta)\wedge1/2} )\rightarrow0$.

The relative speed condition is driven by the biases that arise from
using the discretized observations of $X_t$. The martingale term that
determines the limit behavior of the statistic for a fixed span in
Theorem~\ref{theorem:fixed} is dominated by the empirical process error
$Y_T^{(1)}(u)$ when the time span increases. The leading biases due to
the discretization are two: the drift term $\alpha_t$ and the presence
of ``residual'' jump components in $X_t$ in addition to its leading
stable component at high frequencies. The bias in
$\widehat{\mathcal{L}}_{\beta}(u)$ due to the ``residual'' jump
components is $O_p(|\log\Delta_n|\Delta_n^{1-\beta'/\beta})$. The
higher the activity of the ``residual'' jump components is, the
stronger their effect is on measuring the Laplace transform of
$|\sigma_t|^{\beta}$. Typically, $\beta'$ will be determined from a
Taylor expansion of the \Lvb density of the driving martingale around
zero. In this case $\beta' = \beta-1$ and the bias will be bigger for
the higher levels of activity~$\beta$. The bias due to the drift term
is $O_p(\Delta_n^{2-2/\beta})$, and it becomes bigger the lower the
activity is of the driving martingale. This bias can be significantly
reduced if we make use of $\widetilde{V}_T(X,\Delta_n,\beta,u)$ when
estimating $\mathcal{L}_{\beta}(u)$. Finally, the orders of magnitude
of the above biases can be shown to be optimal by deriving exactly the
bias in the simple case (covered by our Assumption~\ref{assa}) in which $X_t$ is
\Lvb and further $Z_t$ is a sum of two independent stable processes
with indexes $\beta$ and $\beta'$.

The relative speed condition here can be compared with the
corresponding one that arises in the problem of maximum likelihood
estimation of Markov jump-diffusions; see, for example,~\cite{SY}. The
general condition in this problem is $T\Delta_n\rightarrow0$ (also
known as the rapidly increasing experimental design); that is, the mesh
of the grid should increase somewhat faster than the time span of the
data. In our problem here we need weaker relative speed condition,
provided we use the stronger Assumption~\ref{assb1} and the deviation of $Z_t$
from a~stable process at high frequencies is not too big; that is,
$\beta'$ is relatively low.

Part b of Theorem~\ref{theorem:clt} makes the limit result in
(\ref{eq:clt_2}) feasible; that is, it provides estimates from the
high-frequency data for the asymptotic variance of the leading term
$Y_T^{(1)}(u)$. The first result in it, that is, the limit in~(\ref{eq:se_0}), is of independent interest. The sample average of the
limit in (\ref{eq:se_0}) essentially identifies the integrated joint
Laplace transform of~$|\sigma_t|^{\beta}$. This is a~natural extension
of our results here for the marginal Laplace transform of~$|\sigma_t|^{\beta}$ and can be used for estimation and testing of the
transitional density specification of the stochastic scale. We do not
pursue this any further here.

Finally, the proof of Theorem~\ref{theorem:clt} implies also that
$\widehat{\mathcal{L}}_{\beta}(u)$ converges to
$\mathcal{L}_{\beta}(u)$ in $L_1(\mathbb{R}_+,\omega)$ where
$\omega(u)$ is a bounded nonnegative-valued weight function with
$\omega(u) = o(u^{-1-\iota})$ when $u\rightarrow\infty$ for arbitrarily
small $\iota>0$. This can be used to invert
$\widehat{\mathcal{L}}_{\beta}(u)$, using regularized kernels as those
of~\cite{Kr03b}, to estimate nonparametrically the density of the
stochastic scale.

\subsection{Long span asymptotics: The case of estimated activity}
The asymptotic results in Theorem~\ref{theorem:clt} rely on the premise
that $\beta$ is known. The realized Laplace transform crucially relies
on $\beta$, as the latter enters not only in its asymptotic limit and
variance but also in its construction. If we put a wrong value of
$\beta$ in the calculation of the realized Laplace transform, then it
is easy to see that $\widehat{\mathcal{L}}_{\beta}(u)$ will converge
either to $1$ or $0$ depending on whether the wrong value is above or
below the true one, respectively.

In this section we provide asymptotic results for the case where the
activity $\beta$ needs to be estimated from the data. Developing an
estimate for $\beta$ from the high-frequency data is relatively easy
(we will give an example at the end of the section). Hence, here we
investigate the effect of estimating~$\beta$ on our asymptotic results
in Theorem~\ref{theorem:clt}.

\begin{theorem}\label{theorem:ts}
Suppose there exists an estimator of $\beta$, denoted with
$\widehat{\beta}$ and Assumptions \textup{\ref{assa}}, \textup{\ref{assb}} and \textup{\ref{assc}} hold.
\begin{longlist}
\item[(a)] If
$\widehat{\beta}-\beta=o_p (\frac{\Delta_n^{\alpha}}{\sqrt
{T}} )$
for some $\alpha>0$, then we have
%
\begin{equation}\label{eq:ts_1}
\sqrt{T} \bigl(\widehat{\mathcal{L}}_{\widehat{\beta}}(u)-\widehat
{\mathcal{L}}_{\beta}(u) \bigr)=o_p \biggl(\frac{1}{\sqrt{T}} \biggr).
\end{equation}
\item[(b)] If $\widehat{\beta}$ uses only information before the
beginning of the sample or an initial part of the sample with a fixed
time-span (i.e., one that does not grow with $T$), and further
$\widehat{\beta}-\beta=O_p(\Delta_n^{\alpha})$ for $\alpha
>1/(2\beta)$,
$\beta'<\beta/2$, $\beta>4/3$ and $T\Delta_n\rightarrow0$, then we
have (locally uniformly in $u$)
%
\begin{equation}\label{eq:ts_2}
\hspace*{25pt}\sqrt{T} \bigl(\widehat{\mathcal{L}}_{\widehat{\beta}}(u)-\widehat
{\mathcal{L}}_{\beta}(u) \bigr)-\frac{\sqrt{T}\log(2u/\Delta_n)
\mathbb{E}(G_{\beta}(u^{1/\beta}|\sigma_t|))}{\beta^2}(\widehat
{\beta}-\beta) \stackrel{\mathbb{P}}{\longrightarrow} 0,\hspace*{-12pt}
\end{equation}
where $G_{\beta}(x)=\beta x^{\beta}e^{-x^{\beta}}$ for $x>0$.
\item
[(c)] Under the conditions of part \textup{(b)}, a consistent estimator
for\break
$\mathbb{E}(G_{\beta}(u\sigma_t))$ is given by
%
%
\begin{eqnarray}\qquad
\widehat{G}_{\beta}&=&\frac{\Delta_n}{T}\sum_{i=1}^{[T/\Delta
_n]} ((2u)^{1/\widehat{\beta}}\Delta_n^{-1/\widehat{\beta
}}\Delta_i^nX )
\sin ((2u)^{1/\widehat{\beta}}\Delta_n^{-1/\widehat{\beta
}}\Delta_i^nX )\nonumber\\ [-8pt]\\ [-8pt]
\qquad& \stackrel{\mathbb{P}}{\longrightarrow}& \mathbb{E}(G_{\beta
}(u^{1/\beta}|\sigma_t|)).\nonumber
\end{eqnarray}
\end{longlist}
\end{theorem}

Unlike the estimation of $\mathcal{L}_{\beta}(u)$, which requires both
$\Delta_n\rightarrow0$ and $T\rightarrow\infty$, the estimation of
$\beta$ can be performed with a fixed time span by only sampling more
frequently. Therefore, typically the error $\widehat{\beta}-\beta$ will
depend only on $\Delta_n$. Thus, in the general case of part (a) of the
theorem, we will need the relative speed condition
$T\Delta_n^{\gamma}\rightarrow0$ for some $\gamma>0$ to guarantee that
the estimation of $\beta$ does not have an asymptotic effect on the
estimation of the Laplace transform of the stochastic scale. By
providing a bit more structure, mainly imposing the restriction that
$\widehat{\beta}$ is estimated by a previous part of the sample or an
initial part of the current sample with a fixed time span, we can
derive the leading component of the introduced error in our estimation.
This is done in part (b) of the theorem, where it is shown that the
latter is a linear function of $\widehat{\beta}-\beta$ (appropriately
scaled). As we mentioned earlier, $\widehat{\beta}$ does not need long
span, just sampling more frequently, that is, $\Delta_n\rightarrow0$.
Therefore, in a practical application one can estimate $\beta$ from a
short period of time at the beginning of the sample and use the
estimated $\widehat{\beta}$ and the rest of the sample (or the whole
sample) to estimate the Laplace transform of the stochastic scale. In
such a case, part (b) allows us to incorporate the asymptotic effect of
the error in estimating $\beta$ into calculation of the standard errors
for $\mathcal{L}_{\beta}(u)$. For this, one needs to note that the
errors in (\ref{eq:clt_1}) and (\ref{eq:ts_2}) in such case are
asymptotically independent.

A more efficient estimator, in the sense of faster rate of convergence,
will mean that the approximation error
$\widehat{\mathcal{L}}_{\widehat{\beta}}(u)-\widehat{\mathcal
{L}}_{\beta}(u)$
will be smaller asymptotically. Finally, the lower bound on $\alpha$ in
part(b) of the above theorem would typically be satisfied when
$\beta'<\beta/2$. We finish this section by providing an example of
$\sqrt{\Delta_n}$-consistent nonparametric estimator of $\beta$ (when
$\beta'<\beta/2$), developed in~\cite{TT09}. The estimation is based on
a ratio of power variations over two time scales for optimally chosen
power. It is formally defined as
%
\begin{equation}\label{eq:bg_est}
\widehat{\beta} =
\frac{\ln (2 )p^*}{\ln (2 )+\ln [\Phi
_T(X,p^*,2\Delta_n) ]-\ln [\Phi_T(X,p^*,\Delta_n) ]},
\end{equation}
where $p^*$ is optimally chosen from a first-step estimation of the
activity, and the power variation $\Phi_T(X,p,\Delta_n)$ is defined as
%
\begin{equation}
\Phi_T(X,p,\Delta_n)=\sum_{i=1}^{[T/\Delta_n]}|\Delta_i^nX|^p.
\end{equation}
It is shown in~\cite{TT09} that $\widehat{\beta}$ in (\ref{eq:bg_est})
is $\sqrt{\Delta_n}$-consistent for $T$ fixed with an associated
feasible central limit theorem also available.

\section{Monte Carlo assessment}\label{sec:mc}
We now examine the properties of the estimators of the Laplace
transform both in the case when activity of the driving martingale is
known or needs to be estimated from the data,
$\widehat{\mathcal{L}}_{\beta}(u)$ and
$\widehat{\mathcal{L}}_{\widehat{\beta}}(u)$, respectively. The Monte
Carlo setup is calibrated for a financial price series. In particular,
we use $1000$ Monte Carlo replications of $1200$ ``days'' worth of
$78$ within-day price increments and this corresponds approximately to
the span and the sampling frequency of our actual data set in the
empirical application. The model used in the Monte Carlo is given
by\looseness=1
%
\begin{equation}\label{eq:model_D}\qquad
dX_t=V_{t}^{1/\beta}\,dL_t,\qquad  dV_t=0.02(1.0-V_t)\,dt+0.05\sqrt{V_t}\,dB_t,
\end{equation}\looseness=0
where $L_t$ is a \Lvb process with characteristic triplet $(0,0,\nu)$
for $\nu(x) = \frac{0.11}{|x|^{1+1.7}}$ or
$\nu(x)=\frac{0.11e^{-0.25|x|}}{|x|^{1+1.7}}$. The first choice of the
\Lvb measure corresponds to that of a stable process with activity
index of $1.7$ while the second one is that of a tempered stable
process with the same value of the activity index. For the second
choice of $\nu(x)$, Assumption~\ref{assa} is satisfied with $\beta' = 0.7$,
which indicates a rather active ``residual'' component in the driving
martingale, in addition to its stable part. Therefore, the second case
represents a~very stringent test for the small sample behavior of the
RLT.

Table~\ref{tb:mc} summarizes the outcome of the Monte Carlo
experiments. The first two columns of the table report the results for
the case when the activity is known and fixed at its true value. In
both cases, the estimate is very accurate and virtually unbiased. The
third column presents the results for the case when the inference is
done with $\beta= 2$ (with $L_t$ being tempered stable) which
corresponds to treating erroneously the process $X_t$ as a
jump-diffusion. As seen from Table~\ref{tb:mc} this results in a rather
nontrivial upward bias. The reason is that in forming the realized
Laplace transform the increments should be inflated by the factor
$\Delta_n^{-1/1.7}$ but they are instead inflated by the much smaller
$\Delta_n^{-1/2}$. Using the under-inflated increments in the
computations induces a very large upward bias in the
estimator.\looseness=1

The last two columns of Table~\ref{tb:mc} summarize the Monte Carlo
results for the case where the index $\beta$ is presumed unknown and
estimated using (\ref{eq:bg_est}) based on the first $252$ ``days'' in
the simulated data set. As to be be expected, the estimator of the
Laplace transform is less accurate than when the activity is known. In
the case when the driving martingale is tempered stable, our measure
becomes slightly biased due to a small bias in the estimate of the
activity level $\beta$. These biases, however, are relatively small
when compared with the standard deviation of the estimator.

%
\begin{table}
\tabcolsep=0pt
\caption{Monte Carlo results} \label{tb:mc}
\begin{tabular*}{\textwidth}{@{\extracolsep{\fill}}lccccc@{}}
\hline
& \textbf{S} & \textbf{TS} & \textbf{TS} & \textbf{S}& \textbf{TS}\\
& \textbf{Fixed at} & \textbf{Fixed at} &
\textbf{Fixed at} & \textbf{Estimated}& \textbf{Estimated}\\
\textbf{Activity}& \textbf{true value} & \textbf{true value} & $\bolds{\beta=2}$\\
\hline
& \multicolumn{5}{c@{}}{$\mathcal{L}_{\beta}(0.10)$} \\
true value & $0.9051$ & $0.9051$ & $0.9051$ & $0.9051$ & $0.9051$ \\
mean & $0.9052$ & $0.9085$ & $0.9390$ & $0.9057$ & $0.9137$ \\
std & $0.0063$ & $0.0065$ & $0.0045$ & $0.0074$ & $0.0072$ \\ [5pt]
& \multicolumn{5}{c@{}}{$\mathcal{L}_{\beta}(0.50)$} \\
true value & $0.6112$ & $0.6112$ & $0.6112$ & $0.6112$ & $0.6112$ \\
mean & $0.6111$ & $0.6159$ & $0.7771$ & $0.6141$ & $0.6434$ \\
std & $0.0208$ & $0.0218$ & $0.0145$ & $0.0292$ & $0.0284$ \\[5pt]
& \multicolumn{5}{c@{}}{$\mathcal{L}_{\beta}(1.25)$} \\
true value & $0.3001$ & $0.3001$ & $0.3001$ & $0.3001$ & $0.3001$ \\
mean & $0.2998$ & $0.3035$ & $0.5776$ & $0.3050$ & $0.3449$ \\
std & $0.0249$ & $0.0259$ & $0.0231$ & $0.0383$ & $0.0393$ \\[5pt]
& \multicolumn{5}{c@{}}{$\mathcal{L}_{\beta}(2.50)$} \\
true value & $0.0980$ & $0.0980$ & $0.0980$ & $0.0974$ & $0.0980$ \\
mean & $0.0977$ & $0.0994$ & $0.3753$ & $0.1024$ & $0.1312$ \\
std & $0.0158$ & $0.0164$ & $0.0265$ & $0.0261$ & $0.0297$ \\[5pt]
& \multicolumn{5}{c@{}}{$\mathcal{L}_{\beta}(3.75)$} \\
true value & $0.0344$ & $0.0344$ & $0.0344$ & $0.0344$ & $0.0344$ \\
mean & $0.0342$ & $0.0350$ & $0.2544$ & $0.0374$ & $0.0536$ \\
std & $0.0084$ & $0.0087$ & $0.0249$ & $0.0142$ & $0.0179$ \\
\hline
\end{tabular*}
\tabnotetext[]{t1}{Note: In all simulated scenarios $T=1200$ and $[1/\Delta_n]=78$.
The mean and the standard deviation (across the Monte Carlo
replications) correspond to the estimator
$\widehat{\mathcal{L}}_{\beta}(u)$ (the first three columns) or
$\widehat{\mathcal{L}}_{\widehat{\beta}}(u)$ (the last two columns).
The estimator $\widehat{\beta}$ is computed using~(\ref{eq:bg_est}) and
the first $252$ ``days'' of the sample. $V_t$ has Gamma marginal law
with corresponding Laplace transform of
$ (1+u*0.05^2/0.04 )^{-0.04/0.05^2}$. The Monte Carlo replica
is $1000$.}
\end{table}

\section{An application to parametric estimation of the stochastic
scale law}\label{sec:param} We apply the preceding theoretical results
to define a criterion for parametric estimation based on contrasting
our nonparametric realized Laplace transform to that of a parametric
model for the stochastic scale (or the time change).

\begin{theorem}\label{theorem:mlt-par}
Suppose the conditions of Theorem~\ref{theorem:clt} are satisfied. Let
the Laplace transform of $|\sigma_t|^{\beta}$ be given by
$\mathcal{L}_{\beta}(u;\theta)$ for some finite-dimensional parameter
vector lying within a compact set $\theta\in\Theta$ with $\theta_0$
denoting the true value and further assume that
$\mathcal{L}_{\beta}(u;\theta)$ is twice continuously-differentiable in
its second argument. If $\Theta^l$ is some local neighborhood of
$\theta_0$, assume
$\sup_{\theta\in\Theta^l} \{|\nabla_{\theta}\mathcal
{L}_{\beta}(u;\theta)|+\nabla_{\theta\theta'}\mathcal{L}_{\beta
}(u;\theta)| \}$
bounded. Suppose for a kernel function with bounded support
$\kappa\dvtx
\mathbb{R}_+\rightarrow\mathbb{R}_+$ we have that
\[
\int_{\mathbb{R}_+} \bigl(\mathcal{L}_{\beta}(u;\theta)-\mathcal
{L}_{\beta}(u;\theta_0) \bigr)^2\kappa(u)\,du>0,\qquad    \theta\neq
\theta_0.
\]

Define the estimator
%
\begin{equation}\label{eq:theta_hat}
\widehat{\theta} =
\operatorname{argmin}\limits_{\theta\in\Theta} \int_{\mathbb{R}_+} \bigl(\widehat
{\mathcal{L}}_{\beta}(u)-\mathcal{L}_{\beta}(u;\theta)
\bigr)^2\widehat{\kappa}(u)\,du,
\end{equation}
where $\widehat{\kappa}$ is a nonnegative estimator of $\kappa$ with
$(u^{2+\iota}\vee
1)\sup_{u\in\mathbb{R}_+}|\widehat{\kappa}(u)-\kappa
(u)| \stackrel{\mathbb{P}}{\longrightarrow}
0$ for some $\iota>0$. Then for $T\rightarrow\infty$ and
$T\Delta_n\rightarrow0$, we have
%
\begin{equation}\label{eq:SE1}\qquad
\sqrt{T} (\widehat{\theta}-\theta_0 ) \stackrel
{\mathcal{L}-s}{\longrightarrow}  \biggl(\int_{\mathbb{R}_+}\nabla
_{\theta}\mathcal{L}_{\beta}(u;\theta)
\nabla_{\theta}\mathcal{L}_{\beta}(u;\theta)'\kappa(u)\,du
\biggr)^{-1}\Xi^{1/2}E',
\end{equation}
where $E'$ is a standard normal vector and
%
\begin{equation}\label{eq:SE2}
\Xi=
\int_{\mathbb{R}_+}\int_{\mathbb{R}_+}\Sigma_{\beta}(u,v)\nabla
_{\theta}\mathcal{L}_{\beta}(u;\theta)\nabla_{\theta}\mathcal
{L}_{\beta}(v;\theta)'\kappa(u)\kappa(v)\,du\,dv,
\end{equation}
for $\Sigma_{\beta}(u,v)$ the variance--covariance of
Theorem~\ref{theorem:clt}.
\end{theorem}

\begin{figure}

\includegraphics{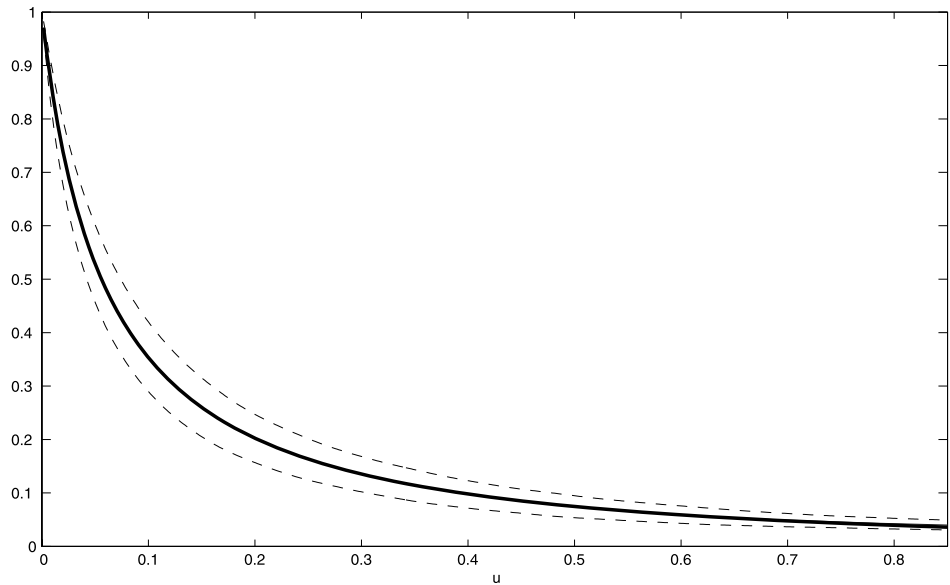}

\caption{The solid line is the fitted parametric
Laplace transform of the marginal distribution of the scale of the log
\VXb index, which essentially plots on top of the estimated realized
Laplace transform (not visible), and the dashed lines show 95 percent
nonparametric confidence intervals about the estimated realized Laplace
transform.} \label{fg:Fit}\vspace*{6pt}
\end{figure}

\begin{remark}
There are two types of pure-jump models
used in practice. First are the time-changed \Lvb processes; see, for
example,~\cite{CGMY03} and~\cite{B11}. As explained in Remark~\ref{rem2}, the
time-change $a_t$ corresponds to $|\sigma_t|^{\beta}$ in
(\ref{eq:intro}). Therefore, $\widehat{\mathcal{L}}_{\beta}(u)$
provides an estimate of the Laplace transform of the time-change which
is modeled directly in parametric settings. The second type of
pure-jump models are the ones specified via \Lv-driven SDE. In this
case we typically model $\sigma_t$ and not $|\sigma_t|^{\beta}$.
Therefore, to apply Theorem~\ref{theorem:mlt-par} in this case one will
need to evaluate $\mathcal{L}_{\beta}(u;\theta)$ via simulation. In
both cases, the use of RLT simplifies the estimation problem
significantly, as it preserves information about the stochastic scale
and, importantly, is robust to any dependence between $\sigma_t$ and
$Z_t$, which, particularly in financial applications, is rather
nontrivial.

The theorem was stated using $\widehat{\mathcal{L}}_{\beta}(u)$, but
obviously the same result will apply if we replace it with
$\widehat{\mathcal{L}}_{\widehat{\beta}}(u)$. By way of illustration,
we apply the theory to the \VXb index computed by the Chicago Board of
Options Exchange; the \VXb is an option-based measure of market
volatility. The data set spans the period from September~22, 2003,
until December 31, 2008, for a total of $1212$ trading days. Within
each day, we use 5-minute records of the \VXb index corresponding to 78
price observations per day.~\cite{TT08} present nonparametric evidence
indicating that the \VXb is a pure-jump \Ito semimartingale.

The underlying pure-jump model we consider for the log \VXb index,
denoted by $v_t$, is
\[
d v_t = \alpha_t\,dt + \int^t_0 \int_{\mathbb{R}} \tilde\mu(dx,dt),
\]
where $\alpha_t$ is the drift term capturing the persistence of $v_t$,
and $ \tilde\mu$ is a random integer-valued measure that has been
compensated by $a_t\, dt \otimes\nu(dx)$ for $a_t$ a stochastic
process capturing time varying intensity. The martingale component of
$v_t$ is a time-changed \Lvb process, as in~\cite{CGMY03}. Recall that
the time-change $a_t$ corresponds to $|\sigma_t|^{\beta}$ in the
general model (\ref{eq:intro}), and our interest here is in making
inferences regarding its marginal distribution.

The parametric specification we use for the marginal distribution of
the time-change is that of a tempered stable subordinator
\cite{Ro04}, which is a self-decomposable distribution, that is, there
is an autoregressive process of order one that generates it
\cite{SATO}. The Laplace transform of the tempered stable is
%
\begin{equation}\label{eq:LP_TS}
\mathcal{L}(u;\theta) =
\cases{
\exp \{c\Gamma(-\alpha) [(\lambda
+u)^{\alpha}-\lambda^{\alpha} ] \},& \mbox{if $\alpha
\in(0,1)$},\vspace*{3pt}\cr
 \displaystyle\biggl( \frac{1}{ 1+u/\lambda} \biggr)^c,& \mbox{if
$\alpha=0$},
}
\end{equation}
where $\theta= (\alpha\,c\,\lambda )$ is the parameter vector,
$\Gamma(-\alpha) = -\frac{1}{\alpha}\Gamma(1-\alpha)$ for
$\alpha\in(0,1)$, $\Gamma$ denotes the standard Gamma function,
$\alpha\in[0,1)$ can be interpreted as the activity index of the
time-change $a_t$, $\,c$ is the scale of the marginal distribution of
$a_t$ and $\lambda$ governs the tail.

To make the estimation feasible, we need an estimate of $\beta$ and to
further specify the kernel $\kappa$ of Theorem~\ref{theorem:mlt-par}.
For $\beta$ we use the estimator defined in equation (\ref{eq:bg_est})
over the first year of the sample, exactly as in the Monte Carlo work;
the point estimate is $\hat\beta= 1.862$ with standard error $0.034$.
We next follow~\cite{Paulson75} in using a Gaussian kernel $\kappa(u) =
\exp(-2u^2/u_{\max}^2)$ where~$u_{\max}$ is defined via $\nabla_u
\mathcal{L}_{\beta}(u_{\max},\theta_0) = -0.05$. The point $u_{\max}$ is
set so that we collect most of the information available in the
empirical Laplace transform. The feasible kernel $\widehat{\kappa}(u)$
is constructed from the infeasible by replacing $u_{\max}$ with a
consistent estimator for it. It is easy to verify that this choice of
the kernels satisfies the conditions of Theorem~\ref{theorem:mlt-par}
above.

\begin{table}[b]
\caption{Estimation results}\label{tb:par_results}
\begin{tabular*}{0.5\textwidth}{@{\extracolsep{\fill}}lcc@{}}
\hline
\textbf{Parameter}& \textbf{Estimate} &\textbf{Standard Error}\\
\hline
$\alpha$ & $0.2651$ & $0.0453$\\
$c$ & $1.2872$ & $0.0469$\\
$\lambda$ & $0.0377$ & $0.0103$\\
\hline
\end{tabular*}
\end{table}

Table~\ref{tb:par_results} shows the parameter estimates and asymptotic
standard errors based on this feasible implementation of
(\ref{eq:theta_hat})--(\ref{eq:SE2}). Interestingly, $\alpha$ is
estimated to be below that associated with the inverse Gaussian
($\alpha=1/2$), while the estimated tail parameter $\lambda$ suggests
relatively moderate dampening, but this parameter appears somewhat
difficult to estimate with high precision, given the time span of our
data set. Figure~\ref{fg:Fit} shows the fit of the model. Specifically,
the heavy solid line is the model-implied Laplace transform evaluated
at the estimated parameters. It plots on top of the (not visible)
realized Laplace transform, (\ref{eq:V}), and thereby passes directly
through the center of the (nonparametric) confidence bands. Overall,
the fit of the tempered stable to the marginal law of the time-change
is quite tight. From this point, one can follow the strategy of
\cite{BNS01} and go further to develop a dynamic model for the
time-change by coupling the fitted marginal law with a specification
for the memory of the process.
\end{remark}

\section{Conclusion}\label{sec:concl}
We derive the asymptotic properties of the realized La\-place transform
for pure-jump processes computed from high-frequency data. The
realized Laplace transform is shown to estimate the Laplace transform
of the stochastic scale of the observed process. The results are
(locally) uniform over the argument of the transform. We can thereby
also derive the asymptotic properties of parameter estimates obtained
by fitting parametric models for the marginal law of the stochastic
scale to the realized Laplace transform. This estimation entails
minimizing a measure of the discrepancy between between the
model-implied and observed transforms.

\section{Proofs}\label{sec:proofs}
Here we give the proof of the main results in the paper:
Lemma~\ref{lema_small_scale} and Theorems~\ref{theorem:fixed} and
\ref{theorem:clt}, with the rest relegated in the supplementary
Appendix. In all the proofs we will denote with $C$ a constant that
does not depend on $T$ and $\Delta_n$, and further it might change from
line to line. We also use the shorthand $\mathbb{E}_{i-1}^n$ for
$\mathbb{E} (\cdot|\mathscr{F}_{(i-1)\Delta_n} )$. We start
with stating some preliminary results, proofs of which are in the
supplement, and which we use in the proofs of the theorems.

\subsection{Preliminary results}\label{subsec:prel}
For a symmetric stable process with \Lvb measure
$\frac{c}{|x|^{\beta+1}}\,dx$ for some $c>0$ and $\beta\in(1,2)$, using
Theorems 14.5 and 14.7 of~\cite{SATO}, we can write its characteristic
function at time $1$ as
\[
\exp \biggl(c\int_0^{\infty}(e^{iur}-1-iur)\frac{dr}{r^{1+\beta
}}+c\int_0^{\infty}(e^{-iur}-1+iur)\frac{dr}{r^{1+\beta}}
\biggr),\qquad  u\in\mathbb{R}.
\]
Then using Lemma 14.11 of~\cite{SATO}, we can simplify the above
expression to
\[
\exp \bigl(2c\Gamma(-\beta)\cos(\beta\pi/2)|u|^{\beta} \bigr),
\]
where $\Gamma(-\beta) = \frac{\Gamma(2-\beta)}{\beta(\beta-1)}$ for
$\beta\in(1,2)$. Therefore, the \Lvb measure of a $\beta$-stable
process, $L_t$, with
$\mathbb{E} (e^{-iuL_t} )=e^{-t|u|^{\beta}/2}$, is
\[
A\times\frac{1}{|x|^{\beta+1}}\,dx,
\]
for $A$ defined in (\ref{eq:A}). Throughout, after appropriately
extending the original probability space, we will use the following
alternative representation of the process $X_t$ (proof of which is
given in the supplement):
\begin{eqnarray*}\label{eq:proof_prel_1}
X_t&=&X_0+\int_0^t\overline{\alpha}_s\,ds+\int_0^t\int_{\mathbb
{R}}\sigma_{s-}x\tilde{\mu}_1(ds,dx)\\
&&{}+\int_0^t\int_{\mathbb
{R}}\sigma_{s-}x\mu_2(ds,dx)
-\int_0^t\int_{\mathbb{R}}\sigma_{s-}x\mu_3(ds,dx)+Y_t,\nonumber
\end{eqnarray*}
where $\mu_1$, $\mu_2$ and $\mu_3$ are homogenous Poisson measures
(the three measures are not mutually independent), with compensators,
respectively, $\nu_1(dx)=\frac{A}{|x|^{\beta+1}}\,dx$,
$\nu_2(dx)=|\nu'(x)|\,dx$ and $\nu_3(dx)=2|\nu'(x)|1 (\nu'(x)<0 )\,dx$ and
$\overline{\alpha}_s =
\alpha_s-\sigma_{s-}\int_{\mathbb{R}}x\nu'(x)\,dx$. Finally, to simplify
notation we will also use the shorthand
$L_t=\int_0^t\int_{\mathbb{R}}x\tilde{\mu}_1(ds,dx)$ and further
for a
symmetric bounded function~$\kappa$ with $\kappa(x)=x$ for $x$ in a
neighborhood of zero, we decompose
%
%
\begin{eqnarray}\label{eq:proof_prel_2}
&&L_t = \overline{L}_t+\widetilde{L}_t,\nonumber\\ [-8pt]\\ [-8pt]
  \eqntext{\displaystyle\overline{L}_t =
\int_{0}^{t}\int_{\mathbb{R}}\bigl(x-\kappa(x)\bigr)\tilde{\mu
}_1(ds,dx),\qquad \widetilde{L}_t
= \int_{0}^{t}\int_{\mathbb{R}}\kappa(x)\tilde{\mu}_1(ds,dx).}
\end{eqnarray}

With this notation we make the following decomposition:
%
%
\begin{eqnarray}\label{eq:proof_prel_3}
 &&V_T(X,\Delta_n,\beta,u)-\int_0^Te^{-u|\sigma_t|^{\beta}}\,dt\nonumber\\ [-8pt]\\ [-8pt]
 &&\qquad=\sum
_{i=1}^{[T/\Delta_n]}\sum_{j=1}^3\xi_{i,u}^{(j)}+\int_{[T/\Delta
_n]\Delta_n}^Te^{-u|\sigma_t|^{\beta}}\,dt,\nonumber
\end{eqnarray}\vspace*{-10pt}
\begin{eqnarray*}
\xi_{i,u}^{(1)} &=&
\Delta_n\cos \bigl((2u)^{1/\beta}\sigma_{(i-1)\Delta_n-}\Delta
_n^{-1/\beta}\Delta_i^nL \bigr)-\int_{(i-1)\Delta_n}^{i\Delta_n}
e^{-u|\sigma_{(i-1)\Delta_n-}|^{\beta}}\,ds,\nonumber\\
\xi_{i,u}^{(2)}
&=&\int_{(i-1)\Delta_n}^{i\Delta_n} \bigl(e^{-u|\sigma_{(i-1)\Delta
_n-}|^{\beta}}-e^{-u|\sigma_{s}|^{\beta}} \bigr)\,ds,\nonumber\\
\xi_{i,u}^{(3)} &=&
\Delta_n \bigl(\cos ((2u)^{1/\beta}\Delta_n^{-1/\beta}\Delta
_i^nX )-\cos \bigl((2u)^{1/\beta}\sigma_{(i-1)\Delta
_n-}\Delta_n^{-1/\beta}\Delta_i^nL \bigr) \bigr).\nonumber
\end{eqnarray*}

Starting with $\xi_{i,u}^{(1)}$, using the self-similarity of the
stable process $L_t$, and the expression for its characteristic
function, we have
%
%
\begin{equation}\label{eq:proof_prel_5}\qquad
\cases{
\mathbb{E}_{i-1}^n \bigl(\cos \bigl((2u)^{1/\beta}\sigma
_{(i-1)\Delta_n-}\Delta_n^{-1/\beta}\Delta_i^nL
 \bigr)-e^{-u|\sigma_{(i-1)\Delta_n-}|^{\beta}} \bigr)=0,\cr
\mathbb{E}_{i-1}^n \bigl(\cos \bigl((2u)^{1/\beta}\sigma
_{(i-1)\Delta_n-}\Delta_n^{-1/\beta}\Delta_i^nL
 \bigr)-e^{-u|\sigma_{(i-1)\Delta_n-}|^{\beta}} \bigr)^2\cr
\qquad=F_{\beta}\bigl(u^{1/\beta}\bigl|\sigma_{(i-1)\Delta_n-}\bigr|\bigr),\cr
\mathbb{E}_{i-1}^n \bigl(\cos \bigl((2u)^{1/\beta}\sigma
_{(i-1)\Delta_n-}\Delta_n^{-1/\beta}\Delta_i^nL
 \bigr)-e^{-u|\sigma_{(i-1)\Delta_n-}|^{\beta}} \bigr)^4\leq
C.
}
\end{equation}
Using first-order Taylor expansion we decompose $\xi_{i,u}^{(2)} =
\sum_{j=1}^3\xi_{i,u}^{(2)}(j)$ where
\begin{eqnarray*}\label{eq:proof_prel_7}
\xi_{i,u}^{(2)}(1)&=&\Upsilon\bigl(\sigma_{(i-1)\Delta_n-},u\bigr)\\
&&{}\times\int
_{(i-1)\Delta_n}^{i\Delta_n}
 \biggl(\int_{(i-1)\Delta_n}^{s}\tilde{\sigma}_u\,dW_u+\int
_{(i-1)\Delta_n}^{s}\int_{\mathbb{R}}\underline{\delta
}(u-,x)\tilde{\underline{\mu}}(du,dx) \biggr)\,ds,\\
\xi_{i,u}^{(2)}(2)&=&\int_{(i-1)\Delta_n}^{i\Delta_n} \bigl(\Upsilon
(\sigma_s^*,u)-\Upsilon\bigl(\sigma_{(i-1)\Delta_n-,u}\bigr) \bigr)\\
 &&{}\times\biggl(\int_{(i-1)\Delta_n}^{s}\tilde{\sigma}_u\,dW_u
+\int_{(i-1)\Delta_n}^{s}\int_{\mathbb{R}^2}\underline{\delta
}(u-,x)\tilde{\underline{\mu}}(du,dx) \biggr)\,ds,\\
\xi_{i,u}^{(2)}(3)&=&\int_{(i-1)\Delta_n}^{i\Delta_n}
\bigl(e^{-u|\widehat{\sigma}_{s}|^{\beta}}-e^{-u|\sigma_{s}|^{\beta
}} \bigr)\,ds,
\end{eqnarray*}
where
$\Upsilon(x,u)=\beta\operatorname{sign}\{x\}u|x|^{\beta-1}e^{-u|x|^{\beta}}$,
$\sigma_s^*$ is a number between $\sigma_{(i-1)\Delta_n-}$ and
$\widehat{\sigma}_s$, and
\begin{eqnarray}\label{eq:proof_prel_8}
\widehat{\sigma}_{s} =
\sigma_{(i-1)\Delta_n-}+\int_{(i-1)\Delta_n}^s\tilde{\sigma
}_u\,dW_u+\int_{(i-1)\Delta_n}^s\int_{\mathbb{R}}\underline{\delta
}(u-,x)\tilde{\underline{\mu}}
(du,dx),\nonumber\\
\eqntext{s\in[(i-1)\Delta_n,i\Delta_n].}
\end{eqnarray}
We derive the following bounds in the supplement for any finite
$\overline{u}>0$
%
\begin{eqnarray}\label{eq:proof_prel_10}
\sum_{i=1}^{[T/\Delta_n]}\mathbb{E}_{i-1}^n \bigl(\xi
_{i,u}^{(2)}(1) \bigr)&=&0,\nonumber\\ [-8pt]\\ [-8pt]
\frac{\Delta_n^{-2}}{T}\mathbb
{E} \Biggl(\sum_{i=1}^{[T/\Delta_n]}\mathbb{E}_{i-1}^n \bigl(\xi
_{i,u}^{(2)}(1) \bigr)^2 \Biggr)&\leq& C, \nonumber\\\label{eq:proof_prel_15}
(T\Delta_n^{\beta/2})^{-1}\sum_{i=1}^{[T/\Delta_n]}\mathbb
{E} \Bigl(\sup_{0\leq u\leq\overline{u}}\bigl|\xi_{i,u}^{(2)}(2)\bigr|
\Bigr)&\leq& C, \\\label{eq:proof_prel_17}
(T\Delta_n)^{-1}\sum_{i=1}^{[T/\Delta_n]}\mathbb{E} \Bigl(\sup
_{0\leq
u\leq\overline{u}}\bigl|\xi_{i,u}^{(2)}(3)\bigr| \Bigr)&\leq&
C.
\end{eqnarray}
Turning to $\xi_{i,u}^{(3)}$, we can first make the following
decomposition [recall the decomposition of $L_t$ in
(\ref{eq:proof_prel_1})]:
%
\begin{equation}\label{eq:chi}
\cos(\chi_1)-\cos(\chi_5) =
\sum_{j=1}^4[\cos(\chi_j)-\cos(\chi_{j+1})],
\end{equation}\vspace*{-15pt}
\begin{eqnarray*}
\chi_1 &=& (2u)^{1/\beta}\Delta_n^{-1/\beta}\Delta_i^nX,\\
 \chi_2 &=&
(2u)^{1/\beta}\Delta_n^{-1/\beta} \biggl(\int_{(i-1)\Delta
_n}^{i\Delta_n}\overline{\alpha}_s\,ds+\int_{(i-1)\Delta_n}^{i\Delta_n}
\sigma_{s-}\,dL_s \biggr),\\
\chi_3 &=&
(2u)^{1/\beta}\Delta_n^{-1/\beta} \biggl(\Delta_n\overline{\alpha
}_{(i-1)\Delta_n}+\int_{(i-1)\Delta_n}^{i\Delta_n}
\sigma_{s-}\,d\widetilde{L}_s \biggr),\\
\chi_4 &=&
(2u)^{1/\beta}\Delta_n^{-1/\beta}\sigma_{(i-1)\Delta_n-}\Delta
_i^n\widetilde{L},\\
 \chi_5
&=&
(2u)^{1/\beta}\Delta_n^{-1/\beta}\sigma_{(i-1)\Delta_n-}\Delta
_i^nL.
\end{eqnarray*}
Then, using the formula for $\cos(x)-\cos(y)$, $x,y\in\mathbb{R}$ for
the first bracketed term on the right-hand side of (\ref{eq:chi}) and a
second-order Taylor expansion for the third one, allows us to write
$\xi_{i,u}^{(3)}=\sum_{j=1}^5\xi_{i,u}^{(3)}(j)$ where
\begin{eqnarray*}\label{eq:proof_prel_19}
\xi_{i,u}^{(3)}(1)&=&-2\Delta_n\\
&&{}\times\sin \biggl(0.5(2u)^{1/\beta}\Delta
_n^{-1/\beta} \biggl(\Delta_i^nX+\int_{(i-1)\Delta_n}^{i\Delta
_n}\overline{\alpha}_s\,ds
+\int_{(i-1)\Delta_n}^{i\Delta_n}\sigma_{s-}\,dL_s \biggr) \biggr)\\
&&{}\times
\sin \biggl(0.5(2u)^{1/\beta}\Delta_n^{-1/\beta} \biggl(\Delta
_i^nX-\int_{(i-1)\Delta_n}^{i\Delta_n}\overline{\alpha}_s\,ds-\int
_{(i-1)\Delta_n}^{i\Delta_n}
\sigma_{s-}\,dL_s \biggr) \biggr),\\\label{eq:proof_prel_20}
\xi_{i,u}^{(3)}(2)&=&-(2u)^{1/\beta}\Delta_n^{2-1/\beta}\sin
\bigl((2u)^{1/\beta}\sigma_{(i-1)\Delta_n-}\Delta_n^{-1/\beta}\Delta
_i^n\widetilde{L} \bigr)
\overline{\alpha}_{(i-1)\Delta_n},\\
\xi_{i,u}^{(3)}(3)&=&-(2u)^{1/\beta}\Delta_n^{1-1/\beta}\sin
\bigl((2u)^{1/\beta}\sigma_{(i-1)\Delta_n-}\Delta_n^{-1/\beta}\Delta
_i^n\widetilde{L} \bigr)\\
&&{}\times\int_{(i-1)\Delta_n}^{i\Delta_n}\bigl(\sigma_{s-}-\sigma_{(i-1)\Delta
_n-}\bigr)\,d\widetilde{L}_s,\\
\xi_{i,u}^{(3)}(4)&=&\Delta_n[\cos(\chi_2)-\cos(\chi_3)] + \Delta
_n[\cos(\chi_4)-\cos(\chi_5)],\\
\xi_{i,u}^{(3)}(5)&=&0.5(2u)^{2/\beta}\Delta_n^{1-2/\beta}\\
&&{}\times\cos
(\tilde{\chi} ) \biggl(\Delta_n\overline{\alpha}_{(i-1)\Delta
_n}+\int_{(i-1)\Delta_n}^{i\Delta_n}
\bigl(\sigma_{s-}-\sigma_{(i-1)\Delta_n-}\bigr)\,d\widetilde{L}_s
\biggr)^2,
\end{eqnarray*}
with $\tilde{\chi}$ denoting some value between $\chi_3$ and $\chi_4$.

We derive the following bounds in the supplement for any finite
$\overline{u}>0$:
%
\begin{eqnarray}\label{eq:proof_prel_22}\qquad
&\hspace*{10pt}\displaystyle (T|\log(\Delta_n)|\Delta_n^{1-\beta'/\beta})^{-1}\sum
_{i=1}^{[T/\Delta_n]}\mathbb{E} \Bigl(\sup_{0\leq u\leq\overline
{u}}\bigl|\xi_{i,u}^{(3)}(1)\bigr| \Bigr)\leq C,\\\label{eq:proof_prel_23}
\qquad&\hspace*{10pt}\displaystyle  \mathbb{E}_{i-1}^n \Bigl(\xi_{i,u}^{(3)}(2) \Bigr)=0,\qquad  \frac
{(\Delta_n^{3-2/\beta})^{-1}}{T}\mathbb{E} \Biggl(\sum
_{i=1}^{[T/\Delta_n]}\mathbb{E}_{i-1}^n \bigl(\xi
_{i,u}^{(3)}(2)\bigr)^2 \Biggr)\leq C,\\\label{eq:proof_prel_25}
\qquad&\hspace*{10pt}\displaystyle  (T\Delta_n^{3/2-1/\beta})^{-1}\sum_{i=1}^{[T/\Delta_n]}\mathbb
{E} \Bigl(\sup_{0\leq u\leq\overline{u}}\bigl|\xi_{i,u}^{(3)}(4)\bigr|
\Bigr)\leq C,\\\label{eq:proof_prel_35}
\qquad&\hspace*{10pt}\displaystyle (T\Delta_n^{2-2/\beta})^{-1}\sum_{i=1}^{[T/\Delta_n]}\mathbb
{E} \Bigl(\sup_{0\leq
u\leq\overline{u}}\bigl|\xi_{i,u}^{(3)}(5)\bigr| \Bigr)\leq
C,\\\label{eq:proof_prel_33}
\qquad&\hspace*{10pt}\displaystyle (T\Delta_n^{1-\iota})^{-1}\mathbb{E} \Biggl(\sup_{0\leq u\leq
\overline{u}} \Biggl|\sum_{i=1}^{[T/\Delta_n]}\mathbb{E}_{i-1}^n\xi
_{i,u}^{(3)}(3) \Biggr| \Biggr)\leq C\\
\eqntext{\mbox{under Assumption~\ref{assb1}},}\\\label{eq:proof_prel_33_a}
\qquad&\hspace*{10pt}\displaystyle  \bigl(T\Delta_n^{1/(\beta\vee\beta^{''}+\iota)}
\bigr)^{-1}\mathbb{E} \Biggl(\sup_{0\leq u\leq\overline{u}} \Biggl|\sum
_{i=1}^{[T/\Delta_n]}\mathbb{E}_{i-1}^n\xi_{i,u}^{(3)}(3)
\Biggr| \Biggr)\leq C\\
\eqntext{\mbox{under Assumption~\ref{assb}},}\\\label{eq:proof_prel_34}
\qquad&\hspace*{10pt}\displaystyle (T\Delta_n^{3-2/\beta})^{-1}\sum_{i=1}^{[T/\Delta_n]}\mathbb
{E} \Bigl(\sup_{0\leq
u\leq\overline{u}}\bigl|\xi_{i,u}^{(3)}(3)\bigr|^2 \Bigr)\leq
C.
\end{eqnarray}
%

\subsection{\texorpdfstring{Proof of Lemma \protect\ref{lema_small_scale}}{Proof of Lemma 1}}
Since $h^{-1/\beta}Z_{ht}$ is a \Lvb process, to prove the convergence
of the sequence, we need to show the convergence of its characteristics
(see, e.g.,~\cite{JS}, Corollary VII.3.6); that is, we need to
establish the following for $h\rightarrow0$:
%
\begin{equation}\label{eq:proof_small_stable_1}
\cases{
\displaystyle h^{1-2/\beta} \int_{\mathbb{R}}\kappa^2(h^{-1/\beta}x)\nu
(x)\,dx\longrightarrow\int_{\mathbb{R}}\kappa^2(x)\frac
{A}{|x|^{\beta+1}}\,dx,\cr
\displaystyle h \int_{\mathbb{R}}g(h^{-1/\beta}x)\nu(x)\,dx\longrightarrow
\int_{\mathbb{R}}g(x)\frac{A}{|x|^{\beta+1}}\,dx,
}
\end{equation}
where $g$ is an arbitrary continuous and bounded function on
$\mathbb{R}$, which is $0$ around~$0$.

The result in (\ref{eq:proof_small_stable_1}) follows by a change of
variable in the integration, and by using the fact that by Assumption~\ref{assa}
we have $|\nu'(x)|<\frac{C}{|x|^{\beta'+1}}$ for $|x|\leq x_0$ where
$x_0$ is fixed and $\beta'<\beta$.

\subsection{\texorpdfstring{Proof of Theorem \protect\ref{theorem:fixed}}{Proof of Theorem 1}}
Part (b) of the theorem holds from the bounds in
(\ref{eq:proof_prel_5})--(\ref{eq:proof_prel_17}) and
(\ref{eq:proof_prel_22})--(\ref{eq:proof_prel_34}), so we are left with
showing part (a). First, we show that for $\Delta_n\rightarrow0$,
$\frac{1}{\sqrt{\Delta_n}}\sum_{i=1}^{[t/\Delta_n]}\xi_{i,u}^{(1)}$
converges stably as a process in~$t$ for the Skorokhod topology to the
process $\int_0^t\sqrt{F_{\beta}(u^{1/\beta}|\sigma_s|)}\,dW_s'$, where~$W_t'$ is a Brownian motion defined on an extension of the original
probability space and is independent from the $\sigma$-field
$\mathcal{F}$. Using the result in (\ref{eq:proof_prel_5}), we get for
every $t>0$,
\[\label{eq:proof_fix_1_1}
\cases{
\displaystyle\frac{1}{\sqrt{\Delta_n}}\sum_{i=1}^{[t/\Delta_n]}\mathbb
{E}_{i-1}^n\bigl(\xi_{i,u}^{(1)}\bigr) \stackrel{\mathbb{P}}{\longrightarrow} 0,\vspace*{3pt}\cr
\displaystyle\frac{1}{\Delta_n}\sum_{i=1}^{[t/\Delta_n]}\mathbb{E}_{i-1}^n\bigl(\xi
_{i,u}^{(1)}\bigr)^2 \stackrel{\mathbb{P}}{\longrightarrow} \int
_0^tF_{\beta}(u^{1/\beta}|\sigma_s|)\,ds,\vspace*{3pt}\cr
\displaystyle\frac{1}{\Delta_n^2}\sum_{i=1}^{[t/\Delta_n]}\mathbb
{E}_{i-1}^n\bigl(\xi_{i,u}^{(1)}\bigr)^4 \stackrel{\mathbb{P}}{\longrightarrow} 0,
}
\]
where for the second convergence above, we made use of
Riemann integrability. Thus to show the stable convergence, given the
above result and upon using Theorem IX.7.28 of~\cite{JS}, we need to
show only
%
\begin{equation}\label{eq:proof_fix_1_2}
\sum_{i=1}^{[t/\Delta_n]}\mathbb{E}_{i-1}^n \bigl(\sqrt{1/\Delta
_n}\xi_{i,u}^{(1)}\Delta_i^nM \bigr) \stackrel{\mathbb
{P}}{\longrightarrow} 0\qquad  \forall
t>0,
\end{equation}
where $M$ is a bounded martingale defined on the original probability
space.

When $M$ is discontinuous martingale, we can argue as follows. First,
we can set $M_t^n = M_{[t/\Delta_n]\Delta_n}$ and $N_t^n =
\sum_{i=1}^{[t/\Delta_n]}\sqrt{1/\Delta_n}\xi_{i,u}^{(1)}$ for any $t$.
With this notation we have\vadjust{\goodbreak}
$[M^n,N^n]_t=\sum_{i=1}^{[t/\Delta_n]} (\sqrt{1/\Delta_n}\xi
_{i,u}^{(1)}\Delta_i^nM )$
and $\langle M^n,N^n\rangle_t =
\sum_{i=1}^{[t/\Delta_n]}\mathbb{E}_{i-1}^n (\sqrt{1/\Delta
_n}\xi_{i,u}^{(1)}\Delta_i^nM )$.
We trivially have that $M^n$ converges (for the Skorokhod topology) to
$M$, and furthermore from the results above the limit of $N^n$, which
we denote here with $N$, is a continuous process. Therefore, using
Corollary VI.3.33(b) of~\cite{JS}, we have that $(M^n,N^n)$ is tight. Then,
using the fact that $M$ is a bounded martingale (and hence it has
bounded jumps), we can apply Corollary VI.6.29 of~\cite{JS} and conclude that the
limit of $[M^n,N^n]$ (up to taking a subsequence) is $[M,N]$. However,
since continuous and pure-jump martingales are orthogonal (see, e.g.,
Definition I.4.11 of~\cite{JS}), we conclude that $[M,N] = 0$. Further, the
difference $[M^n,N^n]-\langle M^n,N^n\rangle$ is a martingale, and
using \Ito isometry, the fact that
$\sqrt{1/\Delta_n}\xi_{i,u}^{(1)}\leq C\sqrt{\Delta_n}$, and the
boundedness of $M$, we have
\begin{eqnarray*}
\mathbb{E} ([M^n,N^n]_t-\langle
M^n,N^n\rangle_t )^2&=&\mathbb{E} \biggl(\sum_{s\leq t}(\Delta
M_s^n\Delta N_s^n)^2 \biggr)\\
&\leq& C\Delta_n\mathbb{E} \biggl(\sum
_{s\leq
t}(\Delta M_s^n)^2 \biggr)\leq C\Delta_n.\nonumber
\end{eqnarray*}
Therefore, $[M^n,N^n]-\langle M^n,N^n\rangle$ converges in probability
to zero, and hence so does $\langle M^n,N^n\rangle$.

When $M$ is a continuous martingale, we can write
$\mathbb{E}_{i-1}^n (\xi_{i,u}^{(1)}\Delta_i^nM )=\mathbb
{E}_{i-1}^n (\Delta_i^n
N\Delta_i^nM )$ where now we denote
$N_t=\mathbb{E}(\xi_{i,u}^{(1)}|\mathscr{F}_t)$ for
$t\in[(i-1)\Delta_n,i\Delta_n]$ (which is obviously a martingale with
respect to the filtration $\mathscr{F}_t$). However, note that
$\xi_{i,u}^{(1)}$ is uniquely determined by
$\mathscr{F}_{(i-1)\Delta_n}$ and the homogenous Poisson measure
$\mu_1$. Therefore, $N_t$ remains a martingale for the coarser
filtration $\mathscr{F}^{*}_t=\mathscr{F}_{(i-1)\Delta_n}\cap
\mathscr{F}_t^{\mu_1}$ for $\mathscr{F}_t^{\mu_1}$ being the filtration
generated by the jump measure $\mu_1$. Then using a martingale
representation for the martingale $(N_{t})_{t\geq(i-1)\Delta_n}$ with
respect to the filtration $\mathscr{F}^{*}_t$ (note $\mu_1$ is a
homogenous Poisson measure), Theorem III.4.34 of~\cite{JS}, we can
represent $N_t$ as a~sum of $\mathscr{F}_{(i-1)\Delta_n}$-adapted
variable and an integral, with respect to $\tilde{\mu}_1$. But then
since pure-jump and continuous martingales are orthogonal, we have
$\mathbb{E}_{i-1}^n (\xi_{i,u}^{(1)}\Delta_i^nM )=0$.

This establishes the stable convergence of
$\frac{1}{\sqrt{\Delta_n}}\sum_{i=1}^{[t/\Delta_n]}\xi_{i,u}^{(1)}$.
Next, the bounds for $\xi_{i,u}^{(2)}$ and $\xi_{i,u}^{(3)}$ in
(\ref{eq:proof_prel_10})--(\ref{eq:proof_prel_17}) and
(\ref{eq:proof_prel_22})--(\ref{eq:proof_prel_35}) imply that\break
$\frac{1}{\sqrt{\Delta_n}}\sum_{i=1}^{[T/\Delta_n]}(\xi
_{i,u}^{(2)}+\xi_{i,u}^{(3)})$
is asymptotically negligible for $\Delta_n\rightarrow0$, and $T$
fixed.

\subsection{\texorpdfstring{Proof of Theorem \protect\ref{theorem:clt}}{Proof of Theorem 3}}

\textit{Part} (a). The proof consists of showing
finite-dimensional convergence in $u$ and tightness of the sequence:

(1) \textit{Finite-dimensional convergence}. First, given
Assumption~\ref{assc}, and using a CLT for stationary and ergodic process
(see\vadjust{\goodbreak}
\cite{JS}, Theorem VIII.3.79), we have for a finite-dimensional vector
$\mathbf{u}$,
%
\begin{equation}\label{eq:proof_clt_1}
\sqrt{T} \biggl(\frac{1}{T}\int_0^Te^{-\mathbf{u}|\sigma_t|^{\beta
}}\,dt-\mathcal{L}_{\beta}(\mathbf{u}) \biggr)
 \stackrel{\mathcal{L}}{\longrightarrow} \Psi,
\end{equation}
where $\Psi$ is a zero-mean normal variable with elements of the
variance--covariance matrix given by $\Sigma_{\beta}(u_i,u_j)$.

Next, the results in Section~\ref{subsec:prel} imply for
$T\rightarrow\infty$ and $\Delta_n\rightarrow0$ under the weaker
Assumption~\ref{assb}.
%
%
\begin{eqnarray}\label{eq:proof_clt_2}
\hspace*{35pt}\frac{1}{\sqrt{T}}\sum_{i=1}^{[T/\Delta_n]}\xi_{i,u}^{(1)} &=&
o_p\bigl(\sqrt{\Delta_n}\bigr),\qquad  \frac{1}{\sqrt{T}}\sum_{i=1}^{[T/\Delta
_n]}\xi_{i,u}^{(2)}
= o_p\bigl(\sqrt{T}\Delta_n^{\beta/2}\vee\sqrt{\Delta_n}\bigr),\nonumber\\ [-8pt]\\ [-8pt]
\hspace*{35pt}\frac{1}{\sqrt{T}}\sum_{i=1}^{[T/\Delta_n]}\xi_{i,u}^{(3)} &=&
o_p \bigl(\sqrt{T}\bigl(|\log(\Delta_n)|\Delta_n^{1-\beta'/\beta}\vee
\Delta_n^{(2-2/\beta)\wedge1/2}\bigr)\vee\sqrt{\Delta_n} \bigr),\nonumber
\end{eqnarray}
with the last one replaced with the weaker
\[
\frac{1}{\sqrt{T}}\sum_{i=1}^{[T/\Delta_n]}\xi_{i,u}^{(3)} =
o_p \bigl(\sqrt{T}\bigl(|\log(\Delta_n)|\Delta_n^{1-\beta'/\beta}\vee
\Delta_n^{2-2/\beta}\bigr)\vee\sqrt{\Delta_n} \bigr),
\]
when the stronger
Assumption~\ref{assb1} holds.

(2) \textit{Tightness}. Let's denote for arbitrary $u,v\geq0$,
\[
z_t = \bigl(e^{-u|\sigma_t|^{\beta}}-\mathcal{L}_{\beta}(u)\bigr) -
\bigl(e^{-v|\sigma_t|^{\beta}}-\mathcal{L}_{\beta}(v)\bigr).
\]
Then, using successive conditioning and Lemma VIII.3.102 in~\cite{JS},
together with the boundedness of $z_t$ and Assumption~\ref{assc}, we get
\begin{eqnarray*}\label{eq:proof_clt_3}
\mathbb{E} \Biggl(\frac{1}{\sqrt{T}}\sum_{t=1}^Tz_t\,dt
\Biggr)^2&=&\frac{1}{T}\int_0^T\int_0^T\mathbb{E} (z_tz_s
)\,ds\,dt\\
&\leq& C
|u^{1/p}-v^{1/p}|\frac{1}{T}\int_0^T\int_0^T\mathbb{E}
(|\sigma_{s\wedge
t}|^{\beta/p}\mathbb{E}(z_{s\vee t}|\mathscr{F}_{s\wedge
t}) )\,ds\,dt\\
&\leq&
C|u^{1/p}-v^{1/p}|^{1+\iota}\frac{1}{T}\int_0^T\int_0^T
\bigl(\alpha_{|t-s|}^{\mathrm{mix}} \bigr)^{1/3-\iota}\,dt\,ds\\
&\leq&
C|u^{1/p}-v^{1/p}|^{1+\iota}\int_0^{\infty} (\alpha_s^{\mathrm{mix}} )^{1/3-\iota}\,ds\leq
C|u^{1/p}-v^{1/p}|^{1+\iota},
\end{eqnarray*}
where $\iota>0$ is the constant of Assumption~\ref{assc} and $p>3$. Using
Theorem 12.3 of~\cite{Billingsley}, the above bound implies the
tightness of the sequence
$\frac{1}{\sqrt{T}}\int_0^T(e^{-u|\sigma_t|^{\beta}}-\mathcal
{L}_{\beta}(u))\,dt$,
and from here we have its convergence for the local uniform topology.

Turning now to
$\frac{1}{\sqrt{\Delta_n}}\sum_{i=1}^{[T/\Delta_n]}\xi
_{i,u}^{(1)}$, we
can use the analog of the result in~(\ref{eq:proof_prel_5}) for
$\xi_{i,u}^{(1)}-\xi_{i,v}^{(1)}$, to get
\[
\mathbb{E} \Biggl(\frac{1}{\sqrt{T\Delta_n}}\sum_{i=1}^{[T/\Delta
_n]}\bigl(\xi_{i,u}^{(1)}-\xi_{i,v}^{(1)}\bigr) \Biggr)^2\leq C|u^{1/\beta
}-v^{1/\beta}|^{2},
\]
for some constant $C$. From here, using Theorem 12.3 of
\cite{Billingsley}, we get the tightness of
$\frac{1}{\sqrt{T\Delta_n}}\sum_{i=1}^{[T/\Delta_n]}\xi_{i,u}^{(1)}$.

Similarly, using the analog of (\ref{eq:proof_prel_23}) applied to
$\xi_{i,u}^{(3)}(2)-\xi_{i,v}^{(3)}(2)$, we have
%
\begin{equation}\label{eq:proof_clt_4}\qquad
\frac{\Delta_n^{-(3-2/\beta)}}{T}\mathbb{E} \Biggl(\sum
_{i=1}^{[T/\Delta_n]} \bigl(\xi_{i,u}^{(3)}(2)-\xi
_{i,v}^{(3)}(2) \bigr) \Biggr)^2\leq
C|u^{1/\beta}-v^{1/\beta}|^2.
\end{equation}
This establishes tightness for
$\Delta_n^{-(3/2-1/\beta)}\sum_{i=1}^{[T/\Delta_n]}\xi_{i,u}^{(3)}(2)$.
We can do exactly the same for
$\Delta_n^{-(3/2-1/\beta)}\sum_{i=1}^{[T/\Delta_n]} (\xi
_{i,u}^{(3)}(3)-\mathbb{E}_{i-1}^n(\xi_{i,u}^{(3)}(3)) )$
using the analog of (\ref{eq:proof_prel_34}) applied to
$\xi_{i,u}^{(3)}(3)-\xi_{i,v}^{(3)}(3)-\mathbb{E}_{i-1}^n(\xi
_{i,u}^{(3)}(3)-\xi_{i,v}^{(3)}(3))$.
Next,
%
\begin{equation}\label{eq:proof_clt_5}
\frac{\Delta_n^{-2}}{T}\mathbb{E} \Biggl(\sum_{i=1}^{[T/\Delta
_n]} \bigl(\xi_{i,u}^{(2)}(1)-\xi_{i,v}^{(2)}(1) \bigr)
\Biggr)^2\leq
C (u-v )^2,
\end{equation}
where we used successive conditioning and further made use of the
inequality
\[\label{eq:proof_clt_6}
|\Upsilon(x,u)-\Upsilon(x,v)|\leq C|x|^{\beta-1}|u-v|,
\qquad x\in\mathbb{R}, u,v\geq0,
\]
which follows from applying first-order Taylor expansion of
$\Upsilon(x,u)$ in its second argument and using the fact that the
derivative of $\Upsilon(x,u)$ in its second argument is bounded by
$C|x|^{\beta-1}$. Therefore,
$\sum_{i=1}^{[T/\Delta_n]}\frac{\Delta_n^{-1}}{\sqrt{T}}\xi
_{i,u}^{(2)}(1)$
is tight on the space of continuous functions equipped with the local
uniform topology.

Next, using the results in Section~\ref{subsec:prel}, it is easy to
show that for any finite $\overline{u}>0$, we have
%
%
\begin{equation}\label{eq:proof_clt_7}
\lim_{\Delta_n\downarrow0,
T\uparrow\infty} \mathbb{P} \Biggl(\sup_{0\leq
u\leq\overline{u}} \Biggl|\sum_{i=1}^{[T/\Delta_n]}(T\Delta
_n^{\beta/2})^{-1}\xi_{i,u}^{(2)}(2) \Biggr|>
\epsilon_n \Biggr)=0\qquad  \forall\epsilon_n\uparrow\infty.\hspace*{-35pt}
\end{equation}
The same holds when, in the above, we replace
$(T\Delta_n^{\beta/2})^{-1}\xi_{i,u}^{(2)}(2)$ with either of the
following terms: $(T\Delta_n)^{-1}\xi_{i,u}^{(2)}(3)$,
$(T|\log(\Delta_n)|\Delta_n^{1-\beta'/\beta})^{-1}\xi_{i,u}^{(3)}(1)$,
$(T\Delta_n^{3/2-1/\beta})^{-1}\xi_{i,u}^{(3)}(4)$,
$(T\Delta_n^{2-2/\beta})^{-1}\xi_{i,u}^{(3)}(5)$ as well as
$(T\Delta_n^{^{1-\iota}})^{-1}\mathbb{E}_{i-1}^n(\xi_{i,u}^{(3)}(3))$
under Assumption~\ref{assb1} and
$(T\Delta_n^{1/(\beta\vee\beta^{''}+\iota)})^{-1}\mathbb
{E}_{i-1}^n(\xi_{i,u}^{(3)}(3))$
under the weaker Assumption~\ref{assb}. This implies that those terms are,
uniformly in $u$, bounded in probability.

\textit{Part} (b). First, (\ref{eq:se_0}) follows directly from
Theorem~\ref{theorem:fixed}, so here we only show (\ref{eq:se_1}). If
we denote for $k\geq0$,
\[\label{eq:proof_clt_8}
C_{k,\beta}(u,v) = \frac{1}{T}\sum_{t=k+1}^T\int_{t-1}^t \bigl(
e^{-u|\sigma_s|^{\beta}}-\mathcal{L}_{\beta}(u) \bigr)\,ds\int_{t-k-1}^{t-k}
 \bigl(e^{-v|\sigma_s|^{\beta}}-\mathcal{L}_{\beta}(v) \bigr)\,ds,
\]
then under our assumptions, by standard arguments (see, e.g.,
Proposition~1 in~\cite{Andrews91}), we have
%
\begin{equation}\label{eq:proof_clt_9}\qquad
C_{0,\beta}(u,v)+\sum_{i=1}^{L_T}\omega(i,L_T)\bigl(C_{i,\beta
}(u,v)+C_{i,\beta}(v,u)\bigr) \stackrel{\mathbb{P}}{\longrightarrow
} \Sigma_{\beta}(u,v).
\end{equation}

Therefore, we are left showing
%
\begin{eqnarray}\label{eq:proof_clt_10}
&&\bigl(\widehat{C}_{0,\beta}(u,v)-C_{0,\beta}(u,v)\bigr)\nonumber\\
&&\qquad{}+\sum
_{i=1}^{L_T}\omega(i,L_T)\bigl(\widehat{C}_{i,\beta}(u,v)+\widehat
{C}_{i,\beta}(v,u)\\
&&\hspace*{70pt}\qquad{}-C_{i,\beta}(u,v)-C_{i,\beta}(v,u)\bigr)
 \stackrel{\mathbb{P}}{\longrightarrow} 0.\nonumber
\end{eqnarray}
We note that for arbitrary $1\leq k\leq T$, we have
\[
\Delta_n\sum_{i=[(k-1)/\Delta_n]+1}^{[k/\Delta_n]}\cos
((2u)^{1/\beta}\Delta_n^{-1/\beta}\Delta_i^nX)\leq1  \quad\mbox{and}\quad
\int_{k-1}^ke^{-u|\sigma_s|^{\beta}}\,ds\leq1.
\]
Hence, for $k=0,1,\ldots,L_T$, we have
\begin{eqnarray*}\label{eq:proof_clt_11}
|\widehat{C}_{k,\beta}(u,v)-C_{k,\beta}(u,v) |&\leq&
\frac{1}{T}\sum_{t=1}^T \biggl|\widehat{Z}_{t,\beta}(u)-
\int_{t-1}^te^{-u|\sigma_s|^{\beta}}\,ds \biggr|\\
&&{}  +\frac{1}{T}\sum
_{t=1}^T \biggl|\widehat{Z}_{t,\beta}(v)-
\int_{t-1}^te^{-v|\sigma_s|^{\beta}}\,ds \biggr|\\
&&{}+ |\widehat
{\mathcal{L}}_{\beta}(u)-\mathcal{L}_{\beta}(u) |\\
&&{}+
 |\widehat{\mathcal{L}}_{\beta}(v)-\mathcal{L}_{\beta}(v) |+O \biggl(\frac{k}{T} \biggr).
\end{eqnarray*}
First, using the CLT result in (\ref{eq:proof_clt_1}), and since
$L_T/\sqrt{T}\rightarrow0$, we have
%
\begin{equation}\label{eq:proof_clt_12}
\sum_{i=1}^{L_T}|\omega(i,L_T)| |\widehat{\mathcal{L}}_{\beta
}(u)-\mathcal{L}_{\beta}(u) | \stackrel{\mathbb
{P}}{\longrightarrow} 0\qquad  \forall
u>0.
\end{equation}
Further, using the stationarity of the process $\sigma_t$ and the
bounds on the moments of the terms $\xi_{i,u}^{(j)}$ derived in
Section~\ref{subsec:prel}, we have for every $t\geq1$,
\[\label{eq:proof_clt_13}
\mathbb{E} \biggl|\widehat{Z}_{t,\beta}(u)-\int_{t-1}^te^{-u|\sigma
_s|^{\beta}}\,ds \biggr|
\leq C \bigl(|\log\Delta_n|\Delta_n^{1-\beta'/\beta}\vee\Delta_n^{
(2-2/\beta)\wedge1/2} \bigr).
\]
Therefore, using the relative speed condition between $L_T$ and
$\Delta_n$ in the theorem, we have
%
\begin{equation}\label{eq:proof_clt_14}
\frac{\sum_{i=1}^{L_T}|\omega(i,L_T)|}{T}\sum_{t=1}^T\mathbb
{E} \biggl|\widehat{Z}_{t,\beta}(u)-\int_{t-1}^te^{-u|\sigma
_s|^{\beta}}\,ds \biggr| \longrightarrow 0\qquad \forall
u>0.\hspace*{-35pt}
\end{equation}
(\ref{eq:proof_clt_12}) and (\ref{eq:proof_clt_14}) imply
(\ref{eq:proof_clt_10}), and this, combined with
(\ref{eq:proof_clt_9}), establishes the result in (\ref{eq:se_1}).

\section*{Acknowledgments}
We would like to thank the Editor, an Associate
Editor and two anonymous referees for many constructive
comments which lead to significant improvements.

\begin{supplement}[id=suppA]
\stitle{Supplement to ``Realized Laplace transforms for pure-jump semimartingales''}
\slink[doi]{10.1214/12-AOS1006SUPP} 
\sdatatype{.pdf}
\sfilename{aos1006\_supp.pdf}
\sdescription{This supplement contains proofs of the preliminary results in
Section~\ref{subsec:prel} as well as the proofs of Theorem~\ref{theorem:fixed-2},~\ref{theorem:ts} and~\ref{theorem:mlt-par}.}
\end{supplement}


\printaddresses

\end{document}